\documentclass[12pt,reqno]{amsart}
\pdfoutput=1
\usepackage{latexsym}
\usepackage{color}
\usepackage{graphicx}
\usepackage{amsmath}
\usepackage{amsfonts}
\newtheorem{theorem}{Theorem}[section]
\newtheorem{lemma}[theorem]{Lemma}
\newtheorem{proposition}[theorem]{Proposition}

\theoremstyle{definition}

\theoremstyle{remark}
\newtheorem{remark}[theorem]{Remark}

\DeclareMathOperator{\esssup}{ess\sup}
\begin{document}
\allowdisplaybreaks
\title[Nonlinear wave equations related to the weak null condition]
{Global existence and blow up for systems of nonlinear wave equations related to the weak null condition}
\author[K.\,Hidano and K.\,Yokoyama]{}
\email{hidano@edu.mie-u.ac.jp}
\email{yokoyama@hus.ac.jp}

\renewcommand{\thefootnote}{\fnsymbol{footnote}}
\footnote[0]{2020\textit{ Mathematics Subject Classification}.
 Primary 35L52; Secondary 35L71}



\keywords{Global existence; Blow up; Systems of nonlinear wave equations; 
Weak null condition}
\maketitle

\centerline{\scshape Kunio Hidano}
\smallskip
{\footnotesize
 \centerline{Department of Mathematics, Faculty of Education, Mie University}
   \centerline{1577 Kurima-machiya-cho Tsu, Mie Prefecture 514-8507, Japan}
} 

\medskip

\centerline{\scshape Kazuyoshi Yokoyama}
\smallskip
{\footnotesize
 \centerline{Hokkaido University of Science}
   \centerline{7-Jo 15-4-1 Maeda, Teine, Sapporo, Hokkaido 006-8585, Japan}
}
%
\baselineskip=0.55cm
\begin{abstract}
We discuss how the higher-order term 
$|u|^q$ $(q>1+2/(n-1))$ has nontrivial effects 
in the lifespan of small solutions to the Cauchy problem 
for the system of nonlinear wave equations 
$$
\partial_t^2 u-\Delta u=|v|^p,
\qquad
\partial_t^2 v-\Delta v=|\partial_t u|^{(n+1)/(n-1)}
+|u|^q
$$
in $n\,(\geq 2)$ space dimensions. 
We show the existence of a certain ``critical curve'' 
in the $pq$-plane such that 
for any $(p,q)$ $(p,q>1)$ lying below the curve, 
nonexistence of global solutions occurs, 
whereas for any $(p,q)$ $(p>1+3/(n-1),\,q>1+2/(n-1))$ 
lying exactly on it, 
this system admits a unique global solution for small data. 
When $n=3$, the discussion for the above system 
with $(p,q)=(3,3)$, which lies on the critical curve, 
has relevance to the study on systems satisfying the weak null condition, 
and we obtain a new result of global existence for such systems. 
Moreover, in the particular case of $n=2$ and $p=4$ 
it is observed that 
no matter how large $q$ is, the higher-order term 
$|u|^q$ never becomes negligible 
and it essentially affects the lifespan of small solutions.
\end{abstract}


\section{Introduction}
This paper is concerned with existence or nonexistence 
of global solutions to the Cauchy problem for 
\begin{equation}\label{eq1}
\begin{cases}
\displaystyle{
\partial_t^2 u-\Delta u
=
|v|^p,}
&
\displaystyle{t>0,\,x\in{\mathbb R}^3},\\
\displaystyle{
\partial_t^2 v-\Delta v
=
(\partial_t u)^2
+
|u|^q,}&
\displaystyle{t>0,\,x\in{\mathbb R}^3},\\
\displaystyle{
u(0)=\varepsilon f_1,\,\,\partial_t u(0)=\varepsilon g_1,}
&{}\\
\displaystyle{
v(0)=\varepsilon f_2,\,\,\partial_t v(0)=\varepsilon g_2.}
&{}
\end{cases}
\end{equation}
Moreover, we are also interested in 
existence of global, small solutions 
to that for systems of the form such as 
\begin{equation}\label{202012271206}
\begin{cases}
\displaystyle{
\partial_t^2 u_1-\Delta u_1
=
A^{ij,\alpha\beta}(\partial_\alpha u_i)(\partial_\beta u_j)
+
C_1(u,\partial u),}
&
\displaystyle{t>0,\,x\in{\mathbb R}^3},\\
\displaystyle{
\partial_t^2 u_2-\Delta u_2
=
(\partial_t u_1)^2
+
C_2(u,\partial u),}&
\displaystyle{t>0,\,x\in{\mathbb R}^3},\\
\displaystyle{
u_1(0)=\varepsilon f_1,\,\,\partial_t u_1(0)=\varepsilon g_1,}
&{}\\
\displaystyle{
u_2(0)=\varepsilon f_2,\,\,\partial_t u_2(0)=\varepsilon g_2,}
&{}
\end{cases}
\end{equation}
where $\partial_0=\partial/\partial t$, 
$\partial_k=\partial/\partial x_k$ $(k=1,2,3)$. 
(In fact, 
systems similar to (\ref{eq1}) will be studied 
in general $n\,(\geq 2)$ space dimensions. See (\ref{09041818}) below. 
Moreover, we will actually discuss systems with more general form 
than (\ref{202012271206}). See (\ref{ModelWeakNull}) below. 
Just for the purpose of illustration, 
we consider (\ref{eq1}) and (\ref{202012271206}).)  
Here, and in the following, we use the summation convention, 
that is, if lowered and raised, 
repeated indices of Greek letters 
and Roman letters are summed from 0 to 3 and 1 to 2, respectively. 
$C_i(u,\partial u)$ $(i=1,2)$ are any homogeneous polynomials of 
degree three in $(u_1,\partial u_1, u_2,\partial u_2)$ 
$(\partial u_i=(\partial_t u_i,\nabla u_i),\,i=1,2)$. 
For each $i,j=1,2$, 
the set of the coefficients 
$\{A^{ij,\alpha\beta}:\alpha,\beta=0,1,2,3\}$ 
satisfies the null condition \cite{KlainermanNull86}, 
that is to say, 
for any $X=(X_0,X_1,X_2,X_3)\in{\mathbb R}^{1+3}$ with 
$X_0^2=X_1^2+X_2^2+X_3^2$, 
we have $A^{ij,\alpha\beta}X_\alpha X_\beta=0$. 

As far as the present authors know, 
this paper is the first in the literature 
to discuss the system (\ref{eq1}). 
In order to give the motivation for studying it, 
we start with recalling the origin of the system (\ref{202012271206}), 
reviewing known related results, and raising an open problem. 
Equation (\ref{202012271206}) is well known 
as one of the most typical examples 
satisfying the weak null condition 
that Lindblad and Rodnianski \cite{LR2003} introduced. 
When the cubic nonlinear terms 
$C_1(u,\partial u)$ and $C_2(u,\partial u)$ are absent 
from (\ref{202012271206}) and therefore the system has the form
\begin{equation}\label{202012281658}
\begin{cases}
\displaystyle{
\partial_t^2 u_1-\Delta u_1
=
A^{ij,\alpha\beta}(\partial_\alpha u_i)(\partial_\beta u_j),}
&
\displaystyle{t>0,\,x\in{\mathbb R}^3},\\
\displaystyle{
\partial_t^2 u_2-\Delta u_2
=
(\partial_t u_1)^2,}&
\displaystyle{t>0,\,x\in{\mathbb R}^3},
\end{cases}
\end{equation}  
we know by the work on semilinear hyperbolic systems due to 
Alinhac \cite{Al2006} and Katayama \cite{Katayama2012} that 
global solutions to (\ref{202012281658}) exist 
for small, smooth data. 
(Concerning detailed study on global behaviors of these solutions, 
see \cite{Katayama2012}. See also \cite{KMS}.) 
Also, in the previous paper \cite{HY2017} 
(see also \cite{HZ2019}) 
the present authors studied a quasi-linear system of wave equations 
including (\ref{202012281658}) as an example 
and obtained a priori estimates strong enough 
to show global solutions for small, smooth data. 
One might think that 
cubic terms have no essential effects in the proof 
of global existence of small solutions 
because they are higher order. 
We would like to note, however, that 
higher-order terms such as $|u|^{p_1}$, 
``mixed'' with lower-order terms such as 
$|\partial_t u|^{p_0}$ for some $p_0<p_1$, 
can happen to have nontrivial effects 
in the lifespan of small solutions 
to nonlinear wave equations. 
See, e.g., \cite{Katayama2001}, \cite{HanZhou2014}, 
and \cite{HWY2016}, concerning how 
$|\partial_t u|^{p_0}+|u|^{p_1}$ 
gives unexpected effect on the lifespan. 
(See also Remarks \ref{202001071138} 
and \ref{202102251148} below.) 
If we allow the presence of the cubic terms, 
it causes $(\partial_t u_1)^2+u_1^3$, 
which may have a similar effect, to appear 
on the right-hand side of the second equation 
in (\ref{202012271206}). 
To the best of the present authors' knowledge, 
it has been open to show 
global, small solutions to (\ref{202012271206}) 
in the presence of the cubic terms. 
The system (\ref{eq1}) with $p=q=3$ is expected 
to serve as a simplified model for 
research in this direction. 

Being motivated in this way, 
we set our first purpose of this paper: 
to investigate how the higher-order term 
such as $|u|^q$ $(q>2)$, 
mixed with the lower-order term $(\partial_t u)^2$, 
has effects in the lifespan of small solutions 
to the simplified system (\ref{eq1}). 
It is shown that 
the point $(p,q)=(3,3)$, 
which is naturally associated with (\ref{202012271206}), 
lies on a certain ``critical'' curve in the $pq$-plane. 
Namely, following the argument of 
Ikeda, Sobajima, and Wakasa \cite{ISW}, 
we show that 
if $p,q\,(>1)$ satisfy $pq<2q+3$, 
then nonexistence of global solutions occurs 
for (\ref{eq1}) for some data, 
no matter how small $\varepsilon\,(>0)$ is. 
We also show that if $pq=2q+3$, $p>5/2$, and $q>2$, 
then 
(\ref{eq1}) admits global solutions 
for sufficiently small $\varepsilon$. 
See Theorems \ref{202012251522} and \ref{GETheorem} below. 
In this sense, 
the equation $pq=2q+3$ describes the ``critical'' curve 
in the $pq$-plane. 
We again note that the point $(p,q)=(3,3)$ lies 
on this curve and therefore it belongs to 
the case of global existence for small data. 

Parenthetically, we also remark that 
existence of global, small solutions to (\ref{eq1}) 
with any $(p,q)$ 
$(p>5/2,q>2)$ 
on the critical curve 
forms a sharp contrast to the results 
for the ``$pq$-system'' 
\begin{equation}\label{202012281752}
\partial_t^2 u-\Delta u
=
|v|^p,
\quad
\partial_t^2 v-\Delta v
=
|u|^q,
\quad
t>0,\,\,x\in{\mathbb R}^3
\end{equation}
discussed in \cite{DSGM}, \cite{Deng1999}, 
\cite{KuboOhta}, \cite{AKT}; 
nonexistence of global solutions occurs for some small data 
and any point $(p,q)$ lying on the critical curve 
for (\ref{202012281752}). 
See Figure\,\ref{202102151047} below. 
It is also worth mentioning that 
whenever $p,q>1+\sqrt 2$, 
the system (\ref{202012281752}) admits global solutions 
for small data; on the other hand, 
the famous exponent $1+\sqrt 2$ of John \cite{John1979} 
plays no major role 
in the two results for (\ref{eq1}).

Let us now state the second purpose of this paper. 
It is to prove global existence of small solutions 
to (\ref{202012271206}) featuring the general cubic terms, 
by exploiting some new insights 
that we get from the detailed study of (\ref{eq1}). 
For the aim, we take into account the new viewpoint that 
the system (\ref{202012271206}) is regarded 
as a mixture of the two critical systems: (\ref{202012281658}) 
and 
\begin{equation}\label{202012281750}
\begin{cases}
\displaystyle{
\partial_t^2 u_1-\Delta u_1
=
u_2^3,}
&
\displaystyle{t>0,\,x\in{\mathbb R}^3},\\
\displaystyle{
\partial_t^2 u_2-\Delta u_2
=
(\partial_t u_1)^2
+
u_1^3,}&
\displaystyle{t>0,\,x\in{\mathbb R}^3}.
\end{cases}
\end{equation}
In the opinion of the present authors, 
this is one of the attractive features 
that the system (\ref{202012271206}) offers. 
Since technical key points 
of getting global, small solutions to 
(\ref{202012281658}) or (\ref{202012281750}) 
are helpful in formulating our strategy 
for this second purpose, 
let us summarize them briefly in the next paragraph.

Let us start with (\ref{202012281658}). 
Global existence is due to the fact that 
for each $i,j=1,2$, 
the set of the coefficients 
$\{A^{ij,\alpha\beta}:\alpha,\beta=0,1,2,3\}$ 
satisfies the null condition 
of Klainerman \cite{KlainermanNull86} 
which creates cancellation and provides more decay. 
Without it, 
``almost global existence'' is the most that 
one can expect in general. 
See, e.g., \cite{JK}, \cite{Rammaha1995}, \cite{Deng1999}. 
We can take advantage of the cancellation 
to rewrite the quadratic term 
$A^{ij,\alpha\beta}(\partial_\alpha u_i)(\partial_\beta u_j)$ 
by using the special derivatives 
$T_k=\partial_k+(x_k/|x|)\partial_t$ 
$(k=1,2,3)$. See, e.g., \cite[pp.\,90--91]{Al2010}. 
After this rewriting, 
the energy-type estimate of 
Alinhac \cite[p.\,92]{Al2010} and 
Lindblad-Rodnianski \cite[Corollary 8.2]{LR2005} 
involving a certain space-time $L^2$ weighted estimate 
for the special derivatives has been used commonly 
in \cite{Al2006}, \cite{Katayama2012}, and \cite{HY2017}, 
as one of the crucial ingredients to close the estimate. 
Let us now turn to (\ref{202012281750}). 
It has the form similar to (\ref{eq1}) with 
the critical point $(p,q)=(3,3)$ which 
admits global solutions 
for small data, as shown below. 
The proof of this global existence result 
uses the Li-Zhou estimate concerning the $L^2({\mathbb R}^3)$ norm 
of some fractional derivatives $|D|^s u$. 
It remains valid for (\ref{202012281750}), 
and we enjoy global solutions to the Cauchy problem 
for (\ref{202012281750}) with small data. 
Our strategy for showing global, small solutions to 
the Cauchy problem (\ref{202012271206}) 
is formulated by combining these key techniques 
employed in the proof of global existence for 
(\ref{202012281658}) or (\ref{202012281750}).

The rest of this introduction is devoted 
to setting the notation, followed by explicit statements of 
two theorems concerning blow up and global existence of 
small solutions to the $n$ $(\geq 2)$ dimensional analogue 
of (\ref{eq1}). See (\ref{09041818}) below. 
In the next section, we prove 
nonexistence of global solutions 
for the system (\ref{09041818}). 
At the end of the section, 
we discuss how the higher-order term 
$|u|^q$ $(q>1+2/(n-1))$ combined with 
$|\partial_t u|^{1+2/(n-1)}$ 
affects the lifespan of small solutions to (\ref{09041818}) 
when $n=2,3$. 
Section 3 is concerned with basic commutation relations between 
generalized derivatives, Sobolev-type or trace-type inequalities, 
and a key linear estimate. 
In Section 4, 
by using these, we prove global existence of small solutions 
to (\ref{09041818}) with $n=2,3$. 
In the final section, we explicitly state our global 
existence theorem for (\ref{202012271206}), 
more precisely, for (\ref{ModelWeakNull}). 
Its proof builds upon observations and techniques exhibited 
in Section 4. 
\subsection{Notation}
We employ the notation 
$\langle\tau\rangle:=\sqrt{1+\tau^2}$ for $\tau\in{\mathbb R}$. 
Following Klainerman \cite{KlainermanNull86}, \cite{Kl87}, 
we introduce several partial differential operators as follows: 
$\partial_0=\partial/\partial t$, 
$\partial_j=\partial/\partial {x_j}$, 
$L_j=t\partial_j+x_j\partial_0 $ $(1\leq j\leq n)$, 
$\Omega_{kl}=x_k\partial_l-x_l\partial_k$ $(1\leq k<l\leq n)$, 
$S=t\partial_0+x\cdot\nabla$. 
These operators 
$\partial_0,\dots,\partial_n$, $L_1,\dots,L_n$, 
$\Omega_{12},\dots,\Omega_{1n},\Omega_{23},\dots,
\Omega_{n-1\,n}$ and $S$ are denoted by 
$\Gamma_0,\dots,\Gamma_\nu$ 
$(\nu:=(n^2+3n+2)/2)$ in this order. 
For a multi-index $\alpha=(\alpha_0,\dots,\alpha_\nu)$, 
$\Gamma^{\alpha_0}_0\cdots\Gamma^{\alpha_\nu}_\nu$ is denoted by 
$\Gamma^\alpha$. 
For convenience, we set $\Omega_{kj}:=-\Omega_{jk}$ if $j<k$, 
and $\Omega_{kj}:=0$ if $k=j$. 
Moreover we will use the operator $|D|:=\sqrt{-\Delta}$. 

It is necessary to define the norm for $1\leq p,\,q<\infty$
\begin{equation}
\|v(\cdot)\|_{p,q}
:=
\left(
      \int_0^\infty
      \biggl(
            \int_{S^{n-1}}|v(r\omega)|^qdS_\omega
      \biggr)^{p/q}r^{n-1}dr
\right)^{1/p}
\end{equation}
with an obvious modification for $p=\infty$ 
or $q=\infty$
\begin{align}
&
\|v(\cdot)\|_{\infty,q}
:=
\sup_{r>0}
\biggl(
\int_{S^{n-1}}
|v(r\omega)|^qdS_\omega
\biggr)^{1/q},\\
&
\|v(\cdot)\|_{p,\infty}
:=
\left(
      \int_0^\infty
      \biggl(
            \sup_{\omega\in{S^{n-1}}}
            |v(r\omega)|
      \biggr)^{p}r^{n-1}dr
\right)^{1/p},
\end{align}
where $r=|x|$, $\omega\in S^{n-1}$. 
These types of norms have been effectively used 
for the existence theory of solutions to fully nonlinear wave equations 
in \cite{LY}, \cite{LZ}. 
Let $N$ be a nonnegative integer and $\Psi$ a characteristic function 
of a set of ${\mathbb R}^n$. 
We define the norm
\begin{equation}\label{202012191637}
\|u(t,\cdot)\|_{\Gamma,N,p,q,\Psi}
:=\sum_{|\alpha|\leq N}
\|\Psi(\cdot)\Gamma^\alpha u(t,\cdot)\|_{p,q}.
\end{equation}
For $\Psi\equiv 1$ in (\ref{202012191637}), 
we omit the subscript $\Psi$. 
If $p=q$, then we omit $q$. 
If $N=0$, then we omit both the subscripts $\Gamma$ and $N$. 
In sum, we abbreviate the notation by simply writing 
$\|u(t,\cdot)\|_{\Gamma,N,p,q,\Psi}$ as 
\begin{align*}
&
\|u(t,\cdot)\|_{\Gamma,N,p,q}, \mbox{~when~}\Psi\equiv 1,\\
&
\|u(t,\cdot)\|_{\Gamma,N,p,\Psi}, \mbox{~when~}p=q,\\
&
\|u(t,\cdot)\|_{\Gamma,N,p}, \mbox{~when~}p=q\mbox{~and~}\Psi\equiv 1,\\
&
\|u(t,\cdot)\|_{p,q,\Psi}, \mbox{~when~}N=0,\\
&
\|u(t,\cdot)\|_{p,q}, \mbox{~when~}N=0\mbox{~and~}\Psi\equiv 1,\\
&
\|u(t,\cdot)\|_{p,\Psi}, \mbox{~when~}N=0\mbox{~and~}p=q,\\
&
\|u(t,\cdot)\|_{p}, \mbox{~when~}N=0,\,p=q,\,\mbox{and~}\Psi\equiv 1.
\end{align*}
In this paper, we will employ this notation of the norm 
by choosing $\Psi=\chi_{1,t}$ or $\Psi=\chi_{2,t}$, where 
$\chi_{1,t}$ is the characteristic function of the set 
$\{x\in{\mathbb R}^n:|x|<(1+t)/2\}$ for any fixed $t>0$, 
and $\chi_{2,t}:=1-\chi_{1,t}$. 
In fact, just for simplicity, 
we will simply write $\chi_{1,t}$ and $\chi_{2,t}$ 
as $\chi_1$ and $\chi_2$, respectively. 
\subsection{Two theorems on the $n$ 
$(\geq 2)$ dimensional analogue of (\ref{eq1})}
We consider the Cauchy problem
\begin{align}
\begin{cases}
(\partial_t^2 - \Delta) u(t, x)
=|v(t, x)|^p, & \\
(\partial_t^2 - \Delta) v(t, x)
=
|\partial_t u(t, x)|^{(n+1)/(n-1)}
+
|u(t, x)|^q, & \\  
u(0, x)=\varepsilon f_1(x),\quad 
\partial_t u(0, x)=\varepsilon g_1(x), &  \\
v(0, x)=\varepsilon f_2(x),\quad 
\partial_t v(0, x)=\varepsilon g_2(x), &
\end{cases}
\label{09041818}
\end{align}
for $x \in {\mathbb R}^n$, $t>0$, 
where 
$(f_i, g_i)\in ({\mathcal S}({\mathbb R}^n))^2$ 
$(i=1, 2)$ and $\varepsilon >0$. 
The exponent $(n+1)/(n-1)$, which is equal to $2$ 
for $n=3$, is well known 
as the critical value of $p$ 
for existence or nonexistence of global, small solutions 
to the nonlinear wave equation 
$\partial_t^2 u-\Delta u=|\partial u|^p$. 
See, e.g., \cite{John1981}, \cite{Sideris1983}, \cite{zhou2001}, 
\cite{HWY2012}, and references cited therein. 
We will show that:
\begin{theorem}\label{202012251522}
Let $n\geq 2$. 
If $p>1$, $q>1$, $4q-(n-1)pq+n+3>0$, 
then there exist no global solutions to 
the Cauchy problem {\rm (\ref{09041818})}, 
provided that
\begin{equation}
\begin{aligned}[c]
&
\int_{{\mathbb R}^n} g_i(x)\,dx > 0 \quad (i=1, 2),
\\
&{\rm supp}\, (u, v, \partial_t u, \partial_t v) 
\subset 
\{ (t, x) \in [0, \infty) 
\times {\mathbb R}^n \,:\, |x| \leq t+1\}. 
\end{aligned}
\label{202102121026} 
\end{equation}
\end{theorem}
\begin{figure}[h]
 \centering
 \includegraphics{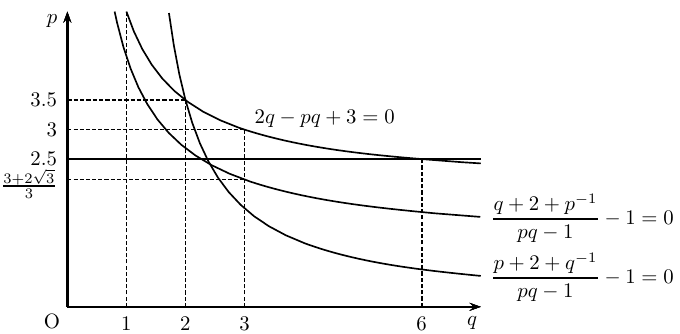}
 \caption{Critical curves in the $pq$-plane 
 related to nonexistence of global solutions when $n=3$. 
$\max \left\{\dfrac{p+2+q^{-1}}{pq-1}, 
\dfrac{q+2+p^{-1}}{pq-1} \right\}-1=0$ is the equation of 
the critical curve for (\ref{202012281752}).
}
\label{202102151047}
\end{figure}
%
In Theorem \ref{202012251522}, 
we say that $(u, v)$ is a solution of (\ref{09041818}) 
for $0 \leq t<T$ if the followings hold:
\begin{enumerate}
\renewcommand{\labelenumi}{(\roman{enumi})}
\item
$\begin{aligned}[t]
&u \in C([0, T);\,H^1({\mathbb R}^n)) \cap C^1([0, T);\,L^2({\mathbb R}^n)) \cap L^q((0,T) \times {\mathbb R}^n),\\
&v \in C([0, T);\,H^1({\mathbb R}^n)) \cap C^1([0, T);\,L^2({\mathbb R}^n)) \cap  L^p((0,T) \times {\mathbb R}^n).
\end{aligned}$
\item
$(u(0), v(0))=(\varepsilon f_1, \varepsilon f_2)$ and $(\partial_t u(0), \partial_t v(0))=(\varepsilon g_1, \varepsilon g_2)$.
\item For every $\varphi \in C_0^{\infty}([0, T) \times {\mathbb R}^n)$,
\begin{align}
&
\begin{aligned}[c]
&\int_0^T\int_{{\mathbb R}^n} |v|^p \varphi\,dxdt \\
&\quad = - \int_{{\mathbb R}^n} \varepsilon g_1(x)\varphi(0, x)\,dx - \int_0^T \int_{{\mathbb R}^n}\left(\partial_t u \partial_t \varphi - \nabla_x u \cdot \nabla_x \varphi\right)\,dxdt,
\end{aligned}
\label{08250926}\\
&
\begin{aligned}[c]
&\int_0^T\int_{{\mathbb R}^n}
\left\{
|\partial_t u|^{(n+1)/(n-1)}+|u|^q
\right\}\varphi\,dxdt\\ 
&\quad = - \int_{{\mathbb R}^n} \varepsilon g_2(x)\varphi(0, x)\,dx - \int_0^T \int_{{\mathbb R}^n}\left(\partial_t v \partial_t \varphi - \nabla_x v \cdot \nabla_x \varphi\right)\,dxdt.
\end{aligned}
\label{08240652} 
\end{align}
\end{enumerate}
Our second theorem on (\ref{09041818}) is 
concerned with global existence of small solutions. 
For a technical reason, 
the space dimension is limited to $n=2$ or $3$. 
We can see by the two theorems that 
the equation $4q-(n-1)pq+n+3=0$ describes the ``critical curve'' 
in the $pq$-plane. 
Set the notation $\sigma(a):=1/2-1/a$ $(a>0)$. 
Using the iteration argument, we prove:
\begin{theorem}\label{GETheorem}
Let $n=2,3$ 
and suppose 
\begin{equation}\label{202102221746}
q
\biggl(p-\frac{4}{n-1}\biggr)
=
1+\frac{4}{n-1},
\quad
p>1+\frac{3}{n-1},
\quad
q>1+\frac{2}{n-1}.
\end{equation}
There exist positive constants 
$C_1$ and $\varepsilon_1$ 
such that 
if $\varepsilon\leq \varepsilon_1$, 
then the Cauchy problem $(\ref{09041818})$ admits a unique global 
solution $(u,v)$ satisfying 
\begin{align}\label{202110031147}
\sum_{|\alpha|\leq 2}
&\esssup\displaylimits_{t>0}
(1+t)^{-\gamma(n,q)}
\|
|D|^{\sigma(q)}\Gamma^\alpha u(t)
\|_{2}
+
\sum_{|\alpha|\leq 2}
\esssup\displaylimits_{t>0}
\|
\partial\Gamma^\alpha u(t)
\|_2
\\
&
+
\sum_{|\alpha|\leq 2}
\esssup\displaylimits_{t>0}
(1+t)^{-1/p}
\|
|D|^{\sigma(p)}\Gamma^\alpha v(t)
\|_{2}\nonumber\\
&
+
\sum_{|\alpha|\leq 2}
\esssup\displaylimits_{t>0}
(1+t)^{-1/(2p)}
\|
|D|^{\sigma(2p)}\Gamma^\alpha v(t)
\|_2
\leq
2C_1\varepsilon_1,
\nonumber
\end{align}
where
\begin{equation*}
\gamma(n,q)
=
\frac{n+1}{2}
\biggl(
\frac{n-1}{n+1}
-
\frac1q
\biggr).
\end{equation*}
\end{theorem}
\begin{remark}
As explained at the end of Section \ref{section202102241600} 
(see Remark \ref{202110081841} below), 
global existence of small solutions to (\ref{09041818}) 
remains true for any exponent pair $(p,q)$ above the critical curve 
(that is, 
the inequality $q(p-4/(n-1))>1+4/(n-1)$ holds 
in place of the equality in (\ref{202102221746})).
\end{remark}
\section{Proof of Theorem \ref{202012251522}}
We apply the test function method as in Ikeda-Sobajima-Wakasa \cite{ISW}. 
Following \cite{ISW}, we take $\eta \in C^{\infty}([0, \infty))$ and $\eta^*$ so that
\begin{align*}
\eta(s)=\begin{cases}
1 & \mbox{for } s<1/2,\\
\mbox{decreasing} & \mbox{for }  1/2 < s < 1,\\
0 & \mbox{for }  s>1,
\end{cases}\quad \mbox{and} \quad
\eta^*(s)=\begin{cases}
0 &\mbox{for }  s<1/2,\\
\eta(s) &\mbox{for }  s\geq 1/2,
\end{cases}
\end{align*}
and define
\begin{align*}
\psi_R(t)=[\eta (R^{-1}t)]^k,\quad
\psi_R^*(t)=[\eta^* (R^{-1}t)]^k
\end{align*}
for $k \geq 2$. We assume that $k \geq \max \{2p', 2q', n+1 \}$. 
Here, $p'$ and $q'$ are 
the H\"older conjugate exponents of 
$p$ and $q$, respectively. 

The following two lemmas are due to Ikeda-Sobajima-Wakasa \cite{ISW} (Lemma 3.2 and Lemma 3.5).

\begin{lemma}\label{202102141027}
Let $1<\bar{p}<\infty$ and $k \geq 2\bar{p}'$. 
Let $w$ be a super-solution of 
$(\partial_t^2 - \Delta) w = H 
\in 
L^1(0, T ; L^1({\mathbb R}^n))$ 
with 
$(w(0), \partial_t w(0)) = (\varepsilon f, \varepsilon g) 
\in (C_0^{\infty}({\mathbb R}^n))^2$ satisfying
\[
{\rm supp}\ w \subset 
\{(t, x) \in [0, \infty) 
\times {\mathbb R}^n : |x| \leq t + r_0\}.
\]
Then 
\[
\varepsilon \int_{{\mathbb R}^n} g(x)\,dx 
+ 
\int_0^T \int_{{\mathbb R}^n} H \psi_R \,dxdt 
\leq 
C R^{-2} \int_0^T \int_{{\mathbb R}^n} |w| [\psi_R^*]^{1/\bar{p}}dxdt
\]
holds for $1 \leq R < T$.
\end{lemma}

\begin{lemma}\label{202102141028}
Let $1<\bar{q}<\infty$ and $k \geq 2\bar{q}'$. 
Let $w$ be a super-solution of 
$(\partial_t^2 - \Delta) w = 0$ 
with 
$(w(0), \partial_t w(0)) = (\varepsilon f, \varepsilon g) 
\in (C_0^{\infty}({\mathbb R}^n))^2$ satisfying
\[
{\rm supp}\ w \subset 
\{(t, x) \in [0, \infty) \times {\mathbb R}^n : |x| \leq t + r_0\},
\quad \int_{{\mathbb R}^n} g(x)\,dx >0.
\]
Then, there exists a constant $\delta>0$ 
independent of $\varepsilon$ and $R$ such that the inequality
\[
\delta 
\left[ 
\varepsilon \int_{{\mathbb R}^n} g(x)\,dx 
\right]^{\bar{q}}  R^{n-(n-1){\bar q}/2} 
\leq 
\int_0^T \int_{{\mathbb R}^n} 
|\partial_t w|^{\bar{q}} \psi_R^*\,dxdt
\]
holds for $1 \leq R < T$.
\end{lemma}

{\it Proof of Theorem $\ref{202012251522}$}. 
We first estimate 
$\int_0^T \int_{{\mathbb R}^n} 
\left\{
|\partial_t u|^{(n+1)/(n-1)}+|u|^q
\right\} \psi_R^*\,dxdt$ from above.
Let $1<R<T$. 
We find from (\ref{08250926}) and Lemma \ref{202102141027} that
\begin{align}
&\int_0^T \int_{{\mathbb R}^n} |v|^p \psi_R\,dxdt
\leq CR^{-2} \int_0^T \int_{{\mathbb R}^n} |u| [\psi_R^*]^{1/q}\,dxdt.
\label{202102121740}
\end{align}
Therefore, we obtain
\begin{equation}
\begin{aligned}[c]
&\int_0^T \int_{{\mathbb R}^n} |v|^p \psi_R\,dxdt\\
&\quad \leq CR^{-2} 
\left( 
\int_0^T \int_{{\mathbb R}^n} |u|^q \psi_R^* \,dxdt 
\right)^{1/q} 
\left( 
\int_{R/2}^R \int_{|x|\leq 1+t}\,dxdt 
\right)^{1/q'} \\
&\quad \leq CR^{-2+(n+1)/q'} 
\left( 
\int_0^T \int_{{\mathbb R}^n} |u|^q \psi_R^*\,dxdt 
\right)^{1/q}.
\end{aligned}
\label{04030838}
\end{equation}
From (\ref{08240652}), we also have
\begin{equation}
\begin{aligned}[c]
&\int_0^T \int_{{\mathbb R}^n} 
\left\{
|\partial_t u|^{(n+1)/(n-1)}
+
|u|^q
\right\} \psi_R\,dxdt\\
&\quad \leq CR^{-2+(n+1)/p'} 
\left( 
\int_0^T \int_{{\mathbb R}^n} |v|^p \psi_R^*\,dxdt 
\right)^{1/p}.
\end{aligned}
\label{09030833}
\end{equation}
Hence, it follows from (\ref{04030838}) 
and (\ref{09030833}) that 
\begin{align}
\left(
\int_0^T \int_{{\mathbb R}^n} 
\left\{|\partial_t u|^{(n+1)/(n-1)}+|u|^q \right\} 
\psi_R^*\,dxdt
\right)^{pq-1} 
&\leq C R^{(n-1)pq-2q-n-1}.
\label{09031224}
\end{align}

Combining (\ref{09031224}) with an estimate of 
$\int_0^T \int_{{\mathbb R}^n} 
\left\{
|\partial_t u|^{(n+1)/(n-1)}+|u|^q
\right\} \psi_R^*\,dxdt$ 
from below, we easily get an upper bound of $R$. 
We use Lemma \ref{202102141028} with 
$\bar{q}=(n+1)/(n-1)$ and obtain
\begin{align}
\delta \varepsilon^{(n+1)/(n-1)} R^{(n-1)/2} 
\leq 
\int_0^T \int_{{\mathbb R}^n} 
|\partial_t u|^{(n+1)/(n-1)} \psi_R^*\,dxdt. 
\label{09031229}
\end{align}
By (\ref{09031224}) and (\ref{09031229}), we obtain 
\[
R \leq C \varepsilon^{-2(n+1)(pq-1)/(n-1)(4q-(n-1)pq+n+3)}.
\] 
Since $R$ is arbitrary as long as $1<R<T$, 
we conclude
\begin{equation}\label{202102251140}
T \leq C \varepsilon^{-2(n+1)(pq-1)/(n-1)(4q-(n-1)pq+n+3)}. 
\end{equation}
We have finished the proof of Theorem \ref{202012251522}. 
\begin{remark}\label{202001071138}
Suppose $n=3$. Let us see 
how $|u|^q+(\partial_t u)^2$ affects 
existence or nonexistence 
of global solutions for small data. 
For the purpose of illustration, 
we focus on $q=3$. Assume (\ref{202102121026}). 
We know by \cite{DSGM} that the Cauchy problem for the system 
\begin{equation*}
\partial_t^2 u-\Delta u=|v|^p,
\quad
\partial_t^2 v-\Delta v=|u|^3,
\quad
t>0,\,\,x\in{\mathbb R}^3
\end{equation*}
admits global solutions 
for small data when $p>(3+2\sqrt 3)/3=2.15\cdots$. 
Moreover, 
we show in Appendix that the Cauchy problem for the system 
\begin{equation*}
\partial_t^2 u-\Delta u=|v|^p,
\quad
\partial_t^2 v-\Delta v=(\partial_t u)^2,
\quad
t>0,\,\,
x\in{\mathbb R}^3
\end{equation*}
with small data admits global solutions when $p>5/2$. 
One might expect that 
global solutions to the system (\ref{09041818}) 
with $p>\max\{5/2,(3+2\sqrt 3)/3\}=5/2$ and $q=3$ 
exist for small data, 
because we are discussing small solutions; 
however, Theorem \ref{202012251522} shows that 
this intuition is false for $5/2<p<3$. 
See Figure\,\ref{202102151047}. 
\end{remark}
\begin{remark}\label{202102251148}
Let us focus on $n=2$ and $p=4$. 
Suitably modifying the proof of 
Proposition \ref{202012171604} and Theorem \ref{202101221717} 
(see (\ref{202101211838}), (\ref{202012211645}) where 
we have $(1+\tau)^{-1}{\hat B}_2(v)^4$ for $n=2$ and $p=4$), 
we easily see that 
if $\varepsilon>0$ is sufficiently small, 
then 
the Cauchy problem 
\begin{align}\label{202102251120}
\begin{cases}
\partial_t^2 u-\Delta u=v^4,&t>0,\,x\in{\mathbb R}^2\\
\partial_t^2 v-\Delta v=|\partial_t u|^3,&t>0,\,x\in{\mathbb R}^2\\  
u(0)=\varepsilon f_1,\quad \partial_t u(0)=\varepsilon g_1, &  \\
v(0)=\varepsilon f_2,\quad \partial_t v(0)=\varepsilon g_2, &
\end{cases}
\end{align}
$(f_i,g_i\in{\mathcal S}({\mathbb R}^2))$ 
admits a unique solution $(u,v)$ defined 
on the interval $[0,T_\varepsilon]$, where 
$T_\varepsilon=\exp(C\varepsilon^{-3})$. 
($C>0$ is a constant independent of $\varepsilon$.) 
We note that for $n=2$ and $p=4$, 
the inequality $4q-(n-1)pq+n+3>0$, 
which appears in Theorem \ref{202012251522}, 
holds for all $q>1$. 
By this, we know that 
there exists a constant $C>1$ depending on $q>1$ 
such that the Cauchy problem 
\begin{align}\label{202102251121}
\begin{cases}
\partial_t^2 u-\Delta u=v^4,&t>0,\,x\in{\mathbb R}^2\\
\partial_t^2 v-\Delta v=|\partial_t u|^3+|u|^q,&t>0,\,x\in{\mathbb R}^2\\  
u(0)=\varepsilon f_1,\quad \partial_t u(0)=\varepsilon g_1, &  \\
v(0)=\varepsilon f_2,\quad \partial_t v(0)=\varepsilon g_2, &
\end{cases}
\end{align}
has no solutions defined 
on the interval 
$[0,C\varepsilon^{-6(4q-1)/5}]$, 
provided that the initial data satisfies (\ref{202102121026}). 
See (\ref{202102251140}). 
By comparing the lifespan for (\ref{202102251120}) 
with that for (\ref{202102251121}), 
we observe that no matter how large $q>1$ is 
and no matter how small $\varepsilon>0$ is, 
the higher-order term $|u|^q$ never becomes negligible 
and it essentially affects the lifespan of small solutions. 
\end{remark}
\section{Commutation and basic inequalities}
Let $[\cdot,\cdot]$ stand for the commutator\,: 
$[A,B]:=AB-BA$. 
\begin{lemma}
The following commutation relations hold 
for $1\leq j<k\leq n$, $l=1,\dots,n$, and $\alpha=0,\dots,n:$ 
\begin{align}
&
[\Omega_{jk},\Box]=0, \quad
[L_l,\Box]=0, \quad
[S,\Box]=-2\Box,\\
&
[S,\Omega_{jk}]=0,\quad
[S,L_l]=0,\quad
[S,\partial_\alpha]=-\partial_\alpha,\\
&
[L_l,\Omega_{jk}]
=
\delta_{lj}L_k-\delta_{lk}L_j.
\end{align}
We also have 
\begin{align}
&
[L_j,L_l]=\Omega_{jl}\quad,1\leq j<l\leq n,\\
&
[L_l,\partial_t]
=
-\partial_l,
\quad
[L_l,\partial_j]
=
-\delta_{lj}\partial_t,\quad l,j=1,\dots,n.
\end{align}
Furthermore, we have for 
$1\leq j<k\leq n,\,1\leq l<m\leq n$
\begin{equation}
[\Omega_{jk},\Omega_{lm}]
=
\delta_{kl}\Omega_{jm}
+
\delta_{km}\Omega_{lj}
+
\delta_{jl}\Omega_{mk}
+
\delta_{jm}\Omega_{kl} 
\end{equation}
and for $1\leq j<k\leq n,\,l=1,\dots,n$
\begin{equation}
[\Omega_{jk},\partial_l]
=
-\delta_{lj}\partial_k
+
\delta_{lk}\partial_j.
\end{equation}
\end{lemma}
The next lemma is concerned with the Klainerman-Sobolev 
type inequalities.
\begin{lemma}Suppose $n\geq 2$. 
$({\rm i})$ Let $2\leq p<q<\infty$, $1/p-1/q\leq 1/n$. Then, 
there exists a positive constant C and the Sobolev-type inequality 
\begin{equation}\label{202012211620}
\|
v(t,\cdot)
\|_{q,\chi_1}
\leq
C(1+t)^{-n(1/p-1/q)}
\|v(t,\cdot)\|_{\Gamma,1,p}
\end{equation}
holds. $({\rm ii})$ Let $2\leq p<\infty$, 
$s\in{\mathbb N}$, and $s>n/p$. 
Then, there exists a positive constant C and the Sobolev-type inequality 
\begin{equation}\label{202012211648}
\|
v(t,\cdot)
\|_{\infty,\chi_1}
\leq
C(1+t)^{-n/p}
\|v(t,\cdot)\|_{\Gamma,s,p}
\end{equation}
holds.
\end{lemma}

It is well known that 
this lemma is an immediate consequence of 
the standard Sobolev embedding and 
\begin{equation}\label{20210124926}
\partial_k
=
\frac{tL_k+\displaystyle{\sum_{j=1}^n}x_j\Omega_{kj}-x_kS}{t^2-r^2},
\quad
k=1,\dots,n,
\end{equation}
together with the standard scaling argument. 

In addition to the Klainerman-Sobolev type inequalities, 
we will make use of the following 
Sobolev-type or trace-type inequalities 
in order to obtain time decay of solutions. 
\begin{lemma}Suppose $n\geq 2$. 
For any $p\in(2,\infty)$ 
there exists a positive constant $C$ and 
the weighted Sobolev-type inequality 
\begin{equation}\label{202012211106}
\||x|^{(n-1)(1/2-1/p)}v\|_{p,2}
\leq
C\||D|^{1/2-1/p}v\|_2
\end{equation}
holds for $v\in C_0^\infty({\mathbb R}^n)$.
\end{lemma}
Since the role and the importance of (\ref{202012211106}) was 
first recognized in \cite{LZ}, 
it has been used in the study of long-time existence 
of small solutions to nonlinear wave equations. 
For the proof, see \cite{LZ} and \cite{DFW}.
\begin{lemma}\mbox{\rm (i)} 
There exists a positive constant $C$ such that 
the trace-type inequality
\begin{equation}\label{20201221932}
r^{1/2}
\|v(r\cdot)\|_{L^4(S^2)}
\leq
C\|\nabla v\|_{L^2({\mathbb R}^3)}
\end{equation}
holds for $v\in C_0^\infty({\mathbb R}^3)$. 
Moreover, for any $s\in (1/2,3/2)$ 
there exists a positive constant $C$ such that 
the trace-type inequality 
\begin{equation}\label{20201221938}
r^{3/2-s}
\|v(r\cdot)\|_{L^p(S^2)}
\leq
C\||D|^s v\|_{L^2({\mathbb R}^3)}
\end{equation}
holds for $v\in C_0^\infty({\mathbb R}^3)$, 
where $2/p=3/2-s$. 
\mbox{\rm (ii)} Let $n=2,3$. 
There exists a positive constant $C$ such that 
the trace-type inequality
\begin{equation}\label{20201221933}
r^{(n-1)/2}\|v(r\cdot)\|_{L^4(S^{n-1})}
\leq
C
\|\nabla v\|_{L^2({\mathbb R}^{n})}^{1/2}
\biggl(
\sum_{|\alpha|\leq 1}
\|\Omega^\alpha v\|_{L^2({\mathbb R}^{n})}
\biggr)^{1/2}
\end{equation}
holds for $v\in C_0^\infty({\mathbb R}^n)$. 
Here, 
$\Omega^\alpha=\Omega_{12}^\alpha$ 
$(\alpha=0,1)$ for $n=2$, 
$\Omega^\alpha
=
\Omega_{12}^{\alpha_1}
\Omega_{13}^{\alpha_2}
\Omega_{23}^{\alpha_3}$ 
$(\alpha=(\alpha_1,\alpha_2,\alpha_3))$ for $n=3$. 
\end{lemma}
For the proof of (\ref{20201221932}) and (\ref{20201221933}) for $n=3$, 
see, e.g., \cite[(3.16) and (3.19)]{Sideris2000}. 
Since we may rely upon the Sobolev embedding 
\begin{equation*}
|v(x)|
\leq
C\sum_{\alpha=0}^1
\|
\Omega^\alpha v(r\cdot)
\|_{L^2(S^1)}, 
\end{equation*}
the proof of (\ref{20201221933}) for $n=2$ is easier than 
that for $n=3$, and we have only to follow the argument 
of \cite[pp.\,864--865]{Sideris2000} and modify a part of it 
suitably. 

Inequality (\ref{20201221938}) is the generalization of 
(\ref{20201221932}). 
It is a direct consequence 
of the Sobolev embedding on $S^2$ and 
the trace-type inequality 
\begin{equation}
r^{3/2-s}
\|v(r\cdot)\|_{L^2(S^2)}
\leq
C\||D|^s\Lambda_\omega^{1/2-s}v\|_{L^2({\mathbb R}^3)}
\end{equation}
obtained in \cite{Hoshiro}, \cite{fangwang}. 
(See \cite{fangwang} as for $\Lambda_\omega$.)

The next lemma is due to Li and Zhou \cite{LZ}, 
and it will play a central role in the proof of 
Theorems \ref{GETheorem} and \ref{GEforLRequation}.
\begin{lemma}Let $n\geq 2$. 
If $\sigma$ satisfies $1/2<1-\sigma<n/2$, 
then the solution $u$ to the inhomogeneous wave equation 
$\partial_t^2 u-\Delta u=F$ in 
${\mathbb R}^n\times (0,\infty)$ 
with data $(f,g)$ at $t=0$ satisfies
\begin{align}\label{202012211613}
&
\||D|^\sigma u(t,\cdot)\|_2
\leq
\||D|^\sigma f\|_2
+
\||D|^{\sigma-1}g\|_2
+
C\int_0^t\|F(\tau,\cdot)\|_{p_1,\chi_1}d\tau\\
&
\hspace{2.3cm}
+
C\int_0^t
\langle\tau\rangle^{-(n/2)+1-\sigma}
\|F(\tau,\cdot)\|_{1,p_2,\chi_2}d\tau.\nonumber
\end{align}
Here $p_1$ and $p_2$ are defined as
\begin{equation}
\frac{1}{p_1}
=
\frac12+\frac{1-\sigma}{n},
\quad
\frac{1}{p_2}
=
\frac12
+
\frac{\frac12-\sigma}{n-1}.
\end{equation}
The functions $\chi_1$ and $\chi_2$ denote 
the characteristic functions of 
$\{x\in{\mathbb R}^n:|x|<(1+\tau)/2\}$ and 
$\{x\in{\mathbb R}^n:|x|>(1+\tau)/2\}$, respectively. 
\end{lemma}
\section{Proof of Theorem \ref{GETheorem}}\label{section202102241600}
In the following discussions, 
we will repeatedly use the notation 
$\sigma(a)=1/2-1/a$ $(a>0)$. 
We consider the Cauchy problem (\ref{09041818}), 
where $f_i$,\,$g_i\in{\mathcal S}({\mathbb R}^n)$ 
$(i=1,2)$ are real-valued functions and 
$\varepsilon>0$ is sufficiently small. 

We introduce two sets of functions
\begin{align}\label{DefinitionUi}
{\mathcal U}_1:=
\{
(u&,v)\in
C(
[0,\infty)
;
{\dot H}^{\sigma(q)}({\mathbb R}^n)
\times{\dot H}^{\sigma(p)}({\mathbb R}^n)
)\,:\,\\
&
\partial_j\Gamma^\alpha u,\,
|D|^{\sigma(q)}\Gamma^\alpha u
\in
C([0,\infty);L^2({\mathbb R}^n)),\,\,
0\leq j\leq n,\,|\alpha|\leq 1,\nonumber\\
&
|D|^{\sigma(p)}\Gamma^\alpha v,\,
|D|^{\sigma(2p)}\Gamma^\alpha v
\in
C([0,\infty);L^2({\mathbb R}^n)),\,\,
|\alpha|\leq 1,
\nonumber\\
&
u(0)=\varepsilon f_1,\,\partial_t u(0)=\varepsilon g_1,\,
v(0)=\varepsilon f_2,\,\partial_t v(0)=\varepsilon g_2\},\nonumber
\end{align}
and
\begin{align}
{\mathcal U}_2:=
\{
(u&,v)\in{\mathcal U}_1\,:\,\\
&
\partial_j\Gamma^\alpha u,\,
|D|^{\sigma(q)}\Gamma^\alpha u
\in
L^\infty((0,\infty);L^2({\mathbb R}^n)),\,\,
0\leq j\leq n,\,|\alpha|\leq 2,\nonumber\\
&
|D|^{\sigma(p)}\Gamma^\alpha v,\,
|D|^{\sigma(2p)}\Gamma^\alpha v
\in
L^\infty((0,\infty);L^2({\mathbb R}^n)),\,\,
|\alpha|\leq 2\}.\nonumber
\end{align}
Also, for $i=1,2$ we set 
$N_i(u,v)=A_i(u)+X_i(u)+
{\tilde B}_i(v)+{\hat B}_i(v)$, where 
\begin{align}
A_1(u)&
=
\sum_{|\alpha|\leq 1}
\sup_{t>0}
(1+t)^{-\gamma(n,q)}
\|
|D|^{\sigma(q)}\Gamma^\alpha u(t)
\|_{2},\label{normA18}\\
\gamma(n,q)
&
=
\frac{n+1}{2}
\biggl(
\frac{n-1}{n+1}
-
\frac1q
\biggr),\\
X_1(u)&
=
\sum_{|\alpha|\leq 1}
\sup_{t>0}
\|
\partial\Gamma^\alpha u(t)
\|_2,\label{normX18}\\
{\tilde B}_1(v)&
=
\sum_{|\alpha|\leq 1}
\sup_{t>0}
(1+t)^{-1/p}
\|
|D|^{\sigma(p)}\Gamma^\alpha v(t)
\|_{2},\label{normB18}\\
{\hat B}_1(v)&
=
\sum_{|\alpha|\leq 1}
\sup_{t>0}
(1+t)^{-1/(2p)}
\|
|D|^{\sigma(2p)}\Gamma^\alpha v(t)
\|_2,\label{hatBnorm}
\end{align}
and $A_2(u)$, $X_2(u)$, ${\tilde B}_2(v)$, 
and ${\hat B}_2(v)$ are defined similarly, with 
$\sum_{|\alpha|\leq 1}$ and $\sup_{t>0}$ 
replaced by 
$\sum_{|\alpha|\leq 2}$ and $\esssup\displaylimits_{t>0}$, 
respectively. 
We mention that the condition 
$q>1+2/(n-1)$ is equivalent to $\gamma(n,q)>0$. 
See (\ref{202102221746}). 
The set ${\mathcal U}_1$ is complete 
with the metric defined as 
$d_1((u,v),({\tilde u},{\tilde v}))
:=
N_1(u-{\tilde u},v-{\tilde v})
$. 
For any positive constant $M$, we set
${\mathcal U}_2(M)
:=
\{(u,v)\in{\mathcal U}_2:\,
N_2(u,v)\leq M\}$. 
By the standard argument, 
we see that if $\varepsilon$ is sufficiently small, 
then it is a non-empty, closed subset of ${\mathcal U}_1$. 
Using the iteration argument, we prove Theorem \ref{GETheorem}. 
\begin{remark}\label{RateRemark2020}
As a test, let us consider the system 
\begin{equation}\label{202012221709}
\Box u=0,
\quad
\Box v=|\partial_t u|^{(n+1)/(n-1)},
\quad
t>0,\,x\in{\mathbb R}^n
\end{equation} 
$(n=2,3)$ with ${\mathcal S}$-data. 
We easily see by the Li-Zhou inequality (\ref{202012211613}) that 
the bound 
$\|\partial\Gamma^\alpha u(t)\|_2=O(1)$ $(|\alpha|\leq 2)$ 
implies 
$\||D|^{\sigma(s)}\Gamma^\alpha v(t)\|_2
=
O(t^{1/s})$, 
for any $s>2$ $(n=2)$, 
$s\geq 2$ $(n=3)$, as $t\to\infty$. 
We regard (\ref{09041818}) as small perturbation 
from (\ref{202012221709}), 
and this is the reason why 
we permit the ${\dot H}^{\sigma(p)}$-norm and 
the ${\dot H}^{\sigma(2p)}$-norm of 
$\Gamma^\alpha v(t)$ 
$(|\alpha|\leq 2)$ 
to grow like $t^{1/p}$ and $t^{1/(2p)}$, 
respectively, as $t\to\infty$. 
See (\ref{normB18}) and (\ref{hatBnorm}) above. 
We note that 
this way of giving a bound on the growth of 
homogeneous Sobolev norms 
is similar to the one in the previous paper \cite{HWY2016}, 
where the global existence of small solutions to 
the Cauchy problem for the single equation 
$\Box u=|\partial_t u|^p+|u|^q$ was discussed 
for some ``critical'' $(p,q)$. 
Regarding it as small perturbation of 
the nonlinear equation $\Box u=|\partial_t u|^p$ 
and giving a certain bound on 
the ${\dot H}^{\sigma(q)}$ norm of $\Gamma^\alpha u(t)$ 
in a similar fashion, 
global existence of small solutions was shown 
even for the ``critical'' $(p,q)$. 
\end{remark}
We carry out the standard iteration argument 
$\Box u_0=\Box v_0=0$, 
$\Box u_m=|v_{m-1}|^p$, 
$\Box v_m=|\partial_t u_{m-1}|^{(n+1)/(n-1)}+|u_{m-1}|^q$ 
$(m=1,2,\dots)$ with data 
$(u_m(0),\partial_t u_m(0))=\varepsilon (f_1,g_1)$ 
and 
$(v_m(0),\partial_t v_m(0))=\varepsilon (f_2,g_2)$. 
Owing to the presence of $|u_{m-1}|^q$, 
one might be concerned that 
the ${\dot H}^{\sigma(p)}$-norm 
or 
the ${\dot H}^{\sigma(2p)}$-norm 
of $\Gamma^\alpha v_m(t)$ 
grows more rapidly than $t^{1/p}$ 
or $t^{1/(2p)}$, respectively. 
(see Remark \ref{RateRemark2020} above). 
In fact, this is not the case; we will observe that 
the ${\dot H}^{\sigma(p)}$-norm 
and the ${\dot H}^{\sigma(2p)}$-norm 
of $\Gamma^\alpha v_m(t)$ 
have the $t^{1/p}$-bound and the $t^{1/(2p)}$-bound 
on the growth, respectively, 
as long as 
the equality $q(p-4/(n-1))=1+4/(n-1)$ holds 
and the ${\dot H}^{\sigma(q)}$ norm 
of $\Gamma^\alpha u_{m-1}(t)$ has the $t^{\gamma(n,q)}$-bound 
on the growth as $t\to\infty$. 

Let us first remark that for $|\alpha|\leq 2$, 
it is possible to express 
$(\Gamma^\alpha u_m,\partial_t\Gamma^\alpha u_m)|_{t=0}$ 
and 
$(\Gamma^\alpha v_m,\partial_t\Gamma^\alpha v_m)|_{t=0}$ 
in terms of $\varepsilon f_i$ 
and $\varepsilon g_i$ $(i=1,2)$, 
and get
\begin{align}\label{202012161715}
\sum_{|\alpha|\leq 2}
\bigl(&
\|(\Gamma^\alpha u_m)(0)\|_{{\dot H}^{\sigma(q)}}
+
\|(\partial_t\Gamma^\alpha u_m)(0)\|_{{\dot H}^{\sigma(q)-1}}
+
\|(\partial\Gamma^\alpha u_m)(0)\|_2\\
&
+
\|(\Gamma^\alpha v_m)(0)\|_{{\dot H}^{\sigma(p)}}
+
\|(\partial_t\Gamma^\alpha v_m)(0)\|_{{\dot H}^{\sigma(p)-1}}\nonumber\\
&
+
\|(\Gamma^\alpha v_m)(0)\|_{{\dot H}^{\sigma(2p)}}
+
\|(\partial_t\Gamma^\alpha v_m)(0)\|_{{\dot H}^{\sigma(2p)-1}}
\bigr)
\leq
C\varepsilon,\nonumber
\end{align}
provided that $\varepsilon$ is sufficiently small. 
Let us see how this bound (\ref{202012161715}) is useful in 
the iteration argument. 
Set $u_{-1}\equiv 0$, $v_{-1}\equiv 0$. 
The crucial point in the proof of Theorem \ref{GETheorem} is to prove the following.
\begin{proposition}\label{202012161730}
 For $m\geq 0$, we have
\begin{align}
A_2(u_m)&
\leq
C\varepsilon
+
C
{\tilde B}_2(v_{m-1})^p,\label{2020June09-921}\\
{\tilde B}_2(v_m)+{\hat B}_2(v_m)
&
\leq
C\varepsilon
+
C
X_2(u_{m-1})^{(n+1)/(n-1)}
+
C
A_2(u_{m-1})^q.\label{2020June09-925}
\end{align}
\end{proposition}
\begin{proposition}\label{202012171604}
For $m\geq 0$, we have
\begin{equation}\label{2020June09-937}
X_2(u_m)
\leq
C\varepsilon
+
C
{\hat B}_2(v_{m-1})^p.
\end{equation}
\end{proposition}
{\it Proof of Proposition $\ref{202012161730}$}. 
Just for simplicity, we drop the subscript $m-1$ in this proof. 
Let us start with the proof of (\ref{2020June09-921}). 
Since we use the Li-Zhou inequality (\ref{202012211613}) 
with $\sigma=\sigma(q)$, 
we need to bound 
$\||v(\tau)|^p\|_{\Gamma,2,p_1,\chi_1}$ and 
$\||v(\tau)|^p\|_{\Gamma,2,1,p_2,\chi_2}$, where 
\begin{equation}
\frac{1}{p_1}
=
\frac12
+
\frac{1-\sigma(q)}{n},\,\,
\frac{1}{p_2}
=
\frac12
+
\frac{\frac12-\sigma(q)}{n-1}
=
\frac12
+
\frac{1}{(n-1)q}.
\end{equation}
Setting $p_3$ and $p_4$ according to
\begin{equation}\label{202110071740}
\frac{1}{p_3}
=
\frac12
-
\frac{\sigma(p)}{n},\,\,
\frac{1}{p_1}
=
\frac{1}{p_3}
+
\frac{1}{p_4},
\end{equation}
we use the H\"older inequality, 
the Klainerman-Sobolev inequality (\ref{202012211620}), 
and the Sobolev embedding 
${\dot H}^{\sigma(p)}({\mathbb R}^n)
\hookrightarrow L^{p_3}({\mathbb R}^n)$, 
to get for $|\alpha|\leq 2$, 
\begin{align}\label{202010301816}
\|&
|v(\tau)|^{p-1}\Gamma^\alpha v(\tau)
\|_{p_1,\chi_1}
\leq
\|
v(\tau)
\|_{p_4(p-1),\chi_1}^{p-1}
\|
\Gamma^\alpha v(\tau)
\|_{p_3}\\
&
\leq
C
(1+\tau)^{-n((p-1)/p_3-1/p_4)}
\|
v(\tau)
\|_{\Gamma,1,p_3}^{p-1}
\|v(\tau)\|_{\Gamma,2,p_3}\nonumber\\
&
\leq
C
(1+\tau)^{-1+\gamma(n,q)}
{\tilde B}_1(v)^{p-1}
{\tilde B}_2(v).\nonumber
\end{align}
Here, we have employed the equality 
in (\ref{202102221746}). 
We can handle 
$\||v(\tau)|^{p-2}(\Gamma^\alpha v(\tau))
(\Gamma^\gamma v(\tau))\|_{p_1,\chi_1}$ with 
$|\alpha|,|\gamma|\leq 1$, similarly. 

Turning our attention to 
$\||v(\tau)|^p\|_{\Gamma,2,1,p_2,\chi_2}$, 
we obtain for $|\alpha|\leq 2$
\begin{align}\label{202010301840}
(&1+\tau)^{-n/2+1-\sigma(q)}
\|
|v(\tau)|^{p-1}\Gamma^\alpha v(\tau)
\|_{1,p_2,\chi_2}\\
&
\leq
(1+\tau)^{-n/2+1-\sigma(q)}
\|
v(\tau)
\|_{p,(n-1)q(p-1),\chi_2}^{p-1}
\|\Gamma^\alpha v(\tau)\|_{p,2,\chi_2}\nonumber\\
&
\leq
C
(1+\tau)^{-n/2+1-\sigma(q)}
\biggl(
\sum_{|\beta|\leq 1}
\|
\Omega^\beta v(\tau)
\|_{p,2,\chi_2}
\biggr)^{p-1}
\|v(\tau)\|_{\Gamma,2,p,2,\chi_2}\nonumber\\
&
\leq
C(1+\tau)^{-1+\gamma(n,q)}
{\tilde B}_1(v)^{p-1}
{\tilde B}_2(v),\nonumber
\end{align}
owing to the equality in (\ref{202102221746}). 
Here, we have employed the Sobolev embedding on $S^{n-1}$ 
and the weighted inequality (\ref{202012211106}). 
In the same way as above, it is possible to deal with 
$\||v(\tau)|^{p-2}(\Gamma^\alpha v(\tau))
(\Gamma^\gamma v(\tau))\|_{1,p_2,\chi_2}$ with 
$|\alpha|,|\gamma|\leq 1$. 

We next show (\ref{2020June09-925}). 
Since we use the Li-Zhou inequality (\ref{202012211613}) 
with $\sigma=\sigma(s)$, $s=p,2p$, 
our task is to bound 
$\||u(\tau)|^q\|_{\Gamma,2,q_1,\chi_1}$, 
$\||u(\tau)|^q\|_{\Gamma,2,1,q_2,\chi_2}$, 
$\||\partial_t u(\tau)|^{(n+1)/(n-1)}\|_{\Gamma,2,q_1,\chi_1}$, 
and 
$\||\partial_t u(\tau)|^{(n+1)/(n-1)}\|_{\Gamma,2,1,q_2,\chi_2}$. 
Here we have set 
\begin{equation}\label{202101071736}
\frac{1}{q_1}
=
\frac12
+
\frac{1-\sigma(s)}n,
\qquad
\frac1{q_2}
=
\frac12
+
\frac{1}{(n-1)s},\,\,s=p,2p.
\end{equation}
Setting $q_3$ and $q_4$ according to
\begin{equation}
\frac{1}{q_3}
=
\frac12
-
\frac{\sigma(q)}{n},
\qquad
\frac{1}{q_1}
=
\frac{1}{q_3}
+
\frac{1}{q_4}
\end{equation}
and repeating the same argument as in (\ref{202010301816}) 
and (\ref{202010301840}), we obtain 
\begin{align}
\|&
|u(\tau)|^q
\|_{\Gamma,2,q_1,\chi_1}
+
(1+\tau)^{-n/2+1-\sigma(s)}
\|
|u(\tau)|^q
\|_{\Gamma,2,1,q_2,\chi_2}\\
&
\leq
C(1+\tau)^{-1+1/s}
A_1(u)^{q-1}A_2(u),\quad
s=p,2p.\nonumber
\end{align}
We turn our attention to the two norms 
of $|\partial_t u|^{(n+1)/(n-1)}$. 
Recalling the definition of $q_1$ 
(see (\ref{202101071736}) above), we set $q_5$ as 
\begin{equation}
\frac{1}{q_5}
=
\frac{1-\sigma(s)}{n},\,\,
\mbox{so that}\,\,\,
\frac{1}{q_1}
=
\frac12
+
\frac{1}{q_5}.
\end{equation}
Using the H\"older inequality and 
the Klainerman-Sobolev inequality (\ref{202012211620}), 
we get for $n=3$
\begin{align}\label{202101211831}
\|
(\partial_t u(\tau))^2
\|_{\Gamma,2,q_1,\chi_1}
&
\leq
C
\biggl(
\sum_{|\alpha|\leq 1}
\|
\partial\Gamma^\alpha u(\tau)
\|_{q_5,\chi_1}
\biggr)
\biggl(
\sum_{|\beta|\leq 2}
\|
\partial\Gamma^\beta u(\tau)
\|_2
\biggr)\\
&
\leq
C
(1+\tau)^{-1+1/s}
X_2(u)^2,\,\,s=p,2p.\nonumber
\end{align}
Similarly, we get for $n=2$ 
\begin{align}\label{202102241613}
\|
|\partial_t u(\tau)|^3
\|_{\Gamma,2,q_1,\chi_1}
&
\leq
C
\biggl(
\sum_{|\alpha|\leq 1}
\|
\partial\Gamma^\alpha u(\tau)
\|_{2q_5,\chi_1}
\biggr)^2
\biggl(
\sum_{|\beta|\leq 2}
\|
\partial\Gamma^\beta u(\tau)
\|_2
\biggr)\\
&
\leq
C
(1+\tau)^{-3/2+1/s}
X_2(u)^3,\,\,s=p,2p.\nonumber
\end{align}
Also, recalling the definition of $q_2$ 
(see (\ref{202101071736}) above), 
we obtain for $n=3$, by the Sobolev embedding on $S^2$
\begin{align}\label{202101211833}
(&1+\tau)^{-3/2+(1-\sigma(s))}
\|
(\partial_t u(\tau))^2
\|_{\Gamma,2,1,q_2,\chi_2}\\
&
\leq
C(1+\tau)^{-1+1/s}
\biggl(
\sum_{|\alpha|\leq 1}
\|
\partial\Gamma^\alpha u(\tau)
\|_{2,2s}
\biggr)
\biggl(
\sum_{|\beta|\leq 2}
\|
\partial\Gamma^\beta u(\tau)
\|_2
\biggr)\nonumber\\
&
\leq
C
(1+\tau)^{-1+1/s}
X_2(u)^2,\quad
s=p,2p,\nonumber
\end{align}
and we also get for $n=2$
\begin{align}\label{202102241616}
(&1+\tau)^{-1+(1-\sigma(s))}
\|
|\partial_t u(\tau)|^3
\|_{\Gamma,2,1,q_2,\chi_2}\\
&
\leq
C(1+\tau)^{-\sigma(s)}
\|\partial_t u(\tau)\|_{\infty,\chi_2}
\biggl(
\sum_{|\alpha|\leq 1}
\|
\partial\Gamma^\alpha u(\tau)
\|_{2,s}
\biggr)
\biggl(
\sum_{|\beta|\leq 2}
\|
\partial\Gamma^\beta u(\tau)
\|_2
\biggr)\nonumber\\
&
\leq
C
(1+\tau)^{-1+1/s}
X_2(u)^3,\quad
s=p,2p,\nonumber
\end{align}
where we have used (\ref{20201221933}) 
together with the Sobolev embedding on $S^1$. 
We have finished the proof of Proposition \ref{202012161730}.

{\it Proof of Proposition $\ref{202012171604}$}. 
Since the proof uses the standard energy estimate, 
we need to deal with 
$\||v(\tau)|^p\|_{\Gamma,2,2}$. 
Via the triangle inequality, we estimate it by 
handling $\||v(\tau)|^p\|_{\Gamma,2,2,\chi_1}$ 
and 
$\||v(\tau)|^p\|_{\Gamma,2,2,\chi_2}$, separately. 
We set $p_5$ and $p_6$ according to
\begin{equation}
\frac{1}{p_5}
=
\frac12
-
\frac{\sigma(2p)}{n},
\quad
\frac{1}{p_6}
=
\frac{\sigma(2p)}{n}
=
\frac{1}{\frac{2np}{p-2}}
+
\frac{1}{2np},
\end{equation}
so that 
$1/2=1/p_5+1/p_6$. 
Using the H\"older inequality, 
the Klainerman-Sobolev inequality (\ref{202012211620}), 
and the Sobolev embedding 
${\dot H}^{\sigma(2p)}({\mathbb R}^n)
\hookrightarrow L^{p_5}({\mathbb R}^n)$, we obtain 
\begin{align}\label{202101211838}
\|&|v(\tau)|^p\|_{\Gamma,2,2,\chi_1}\\
&
\leq
C
\|v(\tau)\|_{2np,\chi_1}^{p-2}
\|v(\tau)\|_{\Gamma,1,2np,\chi_1}
\|v(\tau)\|_{\Gamma,2,p_5}\nonumber\\
&
\leq
C
(1+\tau)^{-n(1/p_5-1/(2np))(p-1)}
\|v(\tau)\|_{\Gamma,2,p_5}^p
\leq
C(1+\tau)^{-(n-1)p/2+n/2}{\hat B}_2(v)^p.\nonumber
\end{align}
Moreover, 
using the H\"older inequality, 
the Sobolev embedding on $S^{n-1}$, 
and the weighted inequality (\ref{202012211106})
we get
\begin{align}\label{202012211645}
\|&
|v(\tau)|^p
\|_{\Gamma,2,2,\chi_2}\\
&
\leq
C\|v(\tau)\|_{2p,\infty,\chi_2}^{p-1}
\|v(\tau)\|_{\Gamma,2,2p,2,\chi_2}
+
C\|v(\tau)\|_{2p,\infty,\chi_2}^{p-2}
\|v(\tau)\|_{\Gamma,1,2p,4,\chi_2}^2
\nonumber\\
&
\leq
C
\|v(\tau)\|_{\Gamma,2,2p,2,\chi_2}^p\nonumber\\
&
\leq
C(1+\tau)^{-(n-1)(p/2-1/2)}
\sum_{|\beta|\leq 2}
\|
|D|^{1/2-1/(2p)}\Gamma^\beta v(\tau)
\|_2^p\nonumber\\
&
\leq
C(1+\tau)^{-(n-1)p/2+n/2}{\hat B}_2(v)^p.\nonumber
\end{align}
Since the condition $p>1+3/(n-1)$ 
(see (\ref{202102221746})) is equivalent to 
$-(n-1)p/2+n/2<-1$, 
we have finished the proof of Proposition \ref{202012171604}.

\vspace{0.5cm}

It immediately follows from Propositions 
\ref{202012161730}--\ref{202012171604} 
that there exist constants $\varepsilon_1>0$ and $C_1>0$ such that 
if $\varepsilon\leq\varepsilon_1$, 
then $(u_m,v_m)\in{\mathcal U}_2(2C_1\varepsilon_1)$ 
for all $m\geq 0$. 
To complete the iteration argument, we need the following:
\begin{proposition}\label{202012211609}
For $m\geq 2$, we get
\begin{align}
N_1&(u_{m+1}-u_m,v_{m+1}-v_m)\\
&
\leq
C(2C_1\varepsilon_1)^{p-1}
\bigl(
{\tilde B}_1(v_m-v_{m-1})
+
{\hat B}_1(v_m-v_{m-1})
\bigr)\nonumber\\
&
+
C
(2C_1\varepsilon_1)^{2/(n-1)}
X_1(u_m-u_{m-1})\nonumber\\
&
+
C(2C_1\varepsilon_1)^{q-1}
A_1(u_m-u_{m-1}).\nonumber
\end{align}
\end{proposition}
It suffices to repeat essentially the same argument as in the proof of 
Propositions \ref{202012161730}--\ref{202012171604}. 
Therefore, we may omit the proof. 

If necessary, we choose $\varepsilon_1$ smaller. 
We then obtain
$
N_1(u_{m+1}-u_m,v_{m+1}-v_m)
\leq
(1/2)^{m-1}
N_1(u_2-u_1,v_2-v_1)$ 
for all $m\geq 2$, 
which means that 
$\{(u_m,v_m)\}$ is a Cauchy sequence in ${\mathcal U}_1$. 
Since ${\mathcal U}_2(2C_1\varepsilon_1)$ is a closed subset of 
${\mathcal U}_1$, 
its limit also lies in ${\mathcal U}_2(2C_1\varepsilon_1)$ 
and it is obviously the unique solution 
to (\ref{09041818}) which we have sought for. 
The proof of Theorem \ref{GETheorem} has been finished.
\begin{figure}[h]
 \centering
 \includegraphics{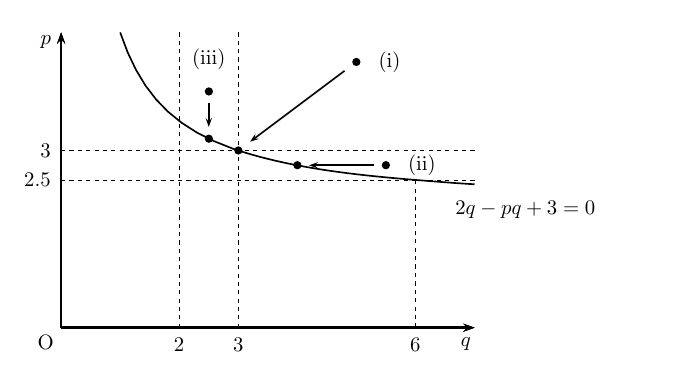}
 \caption{The tail of each arrow represents the point 
 $(p,q)=({\hat p},{\hat q})$. 
 The head of the arrow $({\rm ii})$ represents 
 the point $({\hat p},{\tilde q})$ on the critical curve. 
}
\label{202110121052}
\end{figure}
\begin{remark}\label{202110081841}
Let a pair of the exponent, which we denote by 
$({\hat p},{\hat q})$ in this remark, 
lie above the critical curve, so that 
\begin{equation}
{\hat q}
\biggl({\hat p}-\frac{4}{n-1}\biggr)
>
1+\frac{4}{n-1},
\quad
{\hat p}>1+\frac{3}{n-1},
\quad
{\hat q}>1+\frac{2}{n-1}.
\end{equation}
It is possible to reduce the argument to the critical case 
and obtain global, small solutions to (\ref{09041818}) 
with the exponent $({\hat p},{\hat q})$. 
To verify it, 
we first note the useful $L^\infty$ bound 
for any $(u,v)\in{\mathcal U}_2(M)$. 
Namely, we easily see by using (\ref{202012211648}) and 
(\ref{20201221933}) that 
$\|u(t,\cdot)\|_{\infty}=O(M)$ as $t\to\infty$. 
We also see that 
$\|v(t,\cdot)\|_{\infty}=O(M)$ 
if $p\geq 4$ $(n=2)$, 
$p\geq 3$ $(n=3)$. 
With these $L^\infty$ bounds in mind, for $n=3$, 
we naturally consider the following three cases separately: 
(i) ${\hat p}, {\hat q}>3$, 
(ii) $5/2<{\hat p}\leq 3$, 
(iii) $2<{\hat q}\leq 3$. 
In the case (ii), employing the exponent pair 
$({\hat p}, {\tilde q})$ $({\tilde q}<{\hat q})$ 
on the critical curve (see Figure \ref{202110121052}) 
and choosing $M$ suitably, 
we carry out the iteration argument 
in the set ${\mathcal U}_2(M)$ with 
$(p,q)=({\hat p}, {\tilde q})$, 
with the help of the $L^\infty$ bound for $u$. 
More concretely, this bound enables us to employ 
such a simple inequality as 
$|u|^{{\hat q}-2}|\Gamma_j u||\Gamma_k u|
\leq
CM^{{\hat q}-{\tilde q}}|u|^{{\tilde q}-2}
|\Gamma_j u||\Gamma_k u|$ 
and reduce the argument to the critical case that 
we have just handled. 
In the case (iii), we start with choosing the exponent pair 
$({\tilde p}, {\hat q})$ $({\tilde p}<{\hat p})$ 
on the critical curve. 
Note that when $(u,v)\in{\mathcal U}_2(M)$ with 
$(p,q)=({\tilde p},{\hat q})$, 
we enjoy the $L^\infty$ bound for $v$ 
by virtue of the fact ${\tilde p}\geq 3$. 
It is therefore possible to repeat the above discussion. 
Obviously, the case (i) can be reduced to the case $p=q=3$. 
Finally, we mention that 
the discussion for $n=2$ becomes much simpler, 
because $1+3/(n-1)=4$ for $n=2$. 
We leave the details to the interested reader.
\end{remark}
\section{Model system satisfying the weak null condition}
In this final section, 
we employ the argument in the previous section, 
together with the $L^2$ space-time weighted estimate 
due to Alinhac, Lindblad-Rodnianski for certain special derivatives, 
to discuss one of the most typical systems 
satisfying the weak null condition. 
Let us consider the system of the form
\begin{equation}\label{ModelWeakNull}
\begin{cases}
\displaystyle{
\partial_t^2 u_1-\Delta u_1
=
A_1^{ij,\alpha\beta}
(\partial_\alpha u_i)
(\partial_\beta u_j)
+
C_1(u,\partial u),}
&
\displaystyle{t>0,\,x\in{\mathbb R}^3},\\
\displaystyle{
\partial_t^2 u_2-\Delta u_2
=
A_2^{ij,\alpha\beta}
(\partial_\alpha u_i)
(\partial_\beta u_j)
+
C_2(u,\partial u),}&
\displaystyle{t>0,\,x\in{\mathbb R}^3},\\
\displaystyle{
u_1(0)=\varepsilon f_1,\,\,\partial_t u_1(0)=\varepsilon g_1,}
&{}\\
\displaystyle{
u_2(0)=\varepsilon f_2,\,\,\partial_t u_2(0)=\varepsilon g_2.}
&{}
\end{cases}
\end{equation}
Here, $f_i, g_i\in{\mathcal S}({\mathbb R}^3)$, $i=1,2$. 
Refer to (\ref{202012271206}) 
for the definition of $C_i(u,\partial u)$, $i=1,2$. 
Setting 
${\mathcal N}:=\{X=(X_0,X_1,X_2,X_3)\in{\mathbb R}^4\,:\,
X_0^2=X_1^2+X_2^2+X_3^2\}$, we suppose 
\begin{align}
&
A_1^{ij,\alpha\beta}X_\alpha X_\beta=0
\quad
(i,j=1,2)\,\,\,
\mbox{for all\,$X\in{\mathcal N}$},\label{Aass1}\\
&
A_2^{12,\alpha\beta}X_\alpha X_\beta
=
A_2^{21,\alpha\beta}X_\alpha X_\beta=0\,\,\,
\mbox{for all\,$X\in{\mathcal N}$},\label{Aass2}\\
&
A_2^{22,\alpha\beta}X_\alpha X_\beta=0
\,\,\,
\mbox{for all\,$X\in{\mathcal N}$}.\label{Aass3}
\end{align}
Note that nothing is supposed on $A_2^{11,\alpha\beta}$, 
which in particular permits the presence of $(\partial_t u_1)^2$ 
in the right-hand side of the second equation of (\ref{ModelWeakNull}). 
Under these assumptions, therefore, 
the system (\ref{ModelWeakNull}) is exactly the system 
that Lindblad and Rodnianski introduced 
as one of the models satisfying the weak null condition. 
See \cite[(2), (9), and (13)]{LR2003}. 
To the best of the present authors' knowledge, however, 
it has been open to show the existence of global solutions for 
small and smooth data with compact support. 
In view of Theorem \ref{202012251522}, 
it is obvious that the difficulty comes from 
the coexistence of 
the quadratic term 
$A_2^{11,\alpha\beta}
(\partial_\alpha u_1)
(\partial_\beta u_1)$ in the second equation 
and the cubic terms $u_2^3$ and $u_1^3$ 
in the first and the second equations, respectively. 
Namely, from the point of view of the order of nonlinear terms, 
the system (\ref{ModelWeakNull}) corresponds to 
the critical case 
in the sense that its form is similar to (\ref{eq1}) 
with the critical pair $(p,q)=(3,3)$. 
For the purpose of proving the global existence, 
we naturally employ the norms of the form 
(\ref{normA18})--({\ref{normB18}}) 
with $p=q=3$, 
together with the $L^2$ space-time weighted norm 
(see (\ref{202012171718})--(\ref{DefY(v)}) below) 
to handle the nonlinear terms 
satisfying the null condition (\ref{Aass1})--(\ref{Aass3}). 
(Since we may rely upon the interpolation inequality, 
the norm of the form (\ref{hatBnorm}) is no longer used. 
See (\ref{2021interpolation1}) and (\ref{2021interpolation2}) below.) 
More precisely, 
we suitably modify the definition of ${\mathcal U}_1$ 
(see (\ref{DefinitionUi}) above) and define
\begin{align}
\label{202012181141}{\mathcal V}:=
\{&
(u,v)\in C([0,\infty);H^3\times H^3)\,:\,
N(u,v)<\infty,\\
&
u(0)=\varepsilon f_1,\,\partial_t u(0)=\varepsilon g_1,\,
v(0)=\varepsilon f_2,\,\partial_t v(0)=\varepsilon g_2\},\nonumber
\end{align}
where 
$N(u,v)=W(u)+X(u)+W(v)+Y(v)$, 
\begin{align}
W(u)
=
&\sum_{|\alpha|\leq 2}
\sup_{t>0}
(1+t)^{-1/3}
\|
|D|^{1/6}\Gamma^\alpha u(t)
\|_{2}\\
&
+
\sum_{|\alpha|\leq 3}
\sup_{t>0}
(1+t)^{-1/2-2\delta}
\|\Gamma^\alpha u(t)\|_2\nonumber\\
&
+
\sum_{|\alpha|\leq 2}
\sup_{t>0}
(1+t)^{-1/4}
\||D|^{1/4}\Gamma^\alpha u(t)\|_2,\,\,
\mbox{($W(v)$ is defined similarly),}\nonumber\\
X(u)\label{202012171718}
&=
\sum_{|\alpha|\leq 3}
\biggl(
\sup_{t>0}
\|
\partial\Gamma^\alpha u(t)
\|_2
+
\|
(1+|t-|x||)^{-1/2-\eta}
T\Gamma^\alpha u
\|_{L^2((0,\infty)\times{\mathbb R}^3)}
\biggr),
\\
\label{DefY(v)}
Y(v)&
=
\sum_{|\alpha|\leq 3}
\sup_{t>0}
(1+t)^{-\delta}
\|
\partial\Gamma^\alpha v(t)
\|_2\\
&
+
\sum_{|\alpha|\leq 3}
\sup_{t>0}
(1+t)^{-\delta}
\|
(1+|\tau-|x||)^{-1/2-\eta}
T\Gamma^\alpha v
\|_{L^2((0,t)\times{\mathbb R}^3)}.\nonumber
\end{align}
Here, $\delta$ and $\eta$ are sufficiently small 
positive constants, and 
\begin{equation*}
\|
(1+|t-|x||)^{-1/2-\eta}
T u
\|_{L^2}
:=
\biggl(
\sum_{j=1}^3
\|
(1+|t-|x||)^{-1/2-\eta}
T_j u
\|_{L^2}^2
\biggr)^{1/2}, 
\end{equation*}
$T_j=\partial_j+(x_j/|x|)\partial_t$. 
Using the iteration argument as in the proof 
of Theorem \ref{GETheorem}, 
we prove:
\begin{theorem}\label{GEforLRequation}
Suppose \mbox{$(\ref{Aass1})$--$(\ref{Aass3})$}. 
Then, there exist positive constants 
$C_2$ and $\varepsilon_2$ such that 
if 
$\varepsilon\leq\varepsilon_2$, 
then the Cauchy problem $(\ref{ModelWeakNull})$ 
admits a unique global solution 
$(u_1,u_2)$ satisfying 
$N(u_1,u_2)\leq 2C_2\varepsilon_2$. 
\end{theorem}
Let us compare Theorem \ref{GEforLRequation} with 
\cite[Theorem 1.4]{HY2017} where 
all the cubic and the higher-order terms are absent 
(see \cite[(1.10)]{HY2017}). 
Because of the presence of the general cubic nonlinear terms 
in the system (\ref{ModelWeakNull}), 
the proof of Theorem \ref{GEforLRequation} uses a certain estimate 
for the fractional-order derivatives such as $|D|^{1/6}u$ and $|D|^{1/6}v$, 
which has led us to making the additional assumption (\ref{Aass2}). 
See (\ref{202011181806}). 
On the other hand, 
it immediately follows from \cite[Theorem 1.4]{HY2017} 
that if the cubic nonlinear terms 
$C_1(u,\partial u)$ and $C_2(u,\partial u)$ 
are absent from (\ref{ModelWeakNull}), 
then we have global solutions to (\ref{ModelWeakNull}) 
without assuming (\ref{Aass2}). 
The authors have no idea whether 
global existence result holds without assuming (\ref{Aass2}) 
when the cubic nonlinear terms 
$C_1(u,\partial u)$ and $C_2(u,\partial u)$ 
are present in (\ref{ModelWeakNull}). 

\vspace{0.5cm}

{\it Proof of Theorem $\ref{GEforLRequation}$.} 
In essentially the same way as we have done in the previous section, 
we define $(u_1^{(m)},u_2^{(m)})$ 
$(m=-1,0,\dots)$ inductively. 
To prove Theorem \ref{GEforLRequation}, 
we show:
\begin{proposition}\label{20201030}
For $m=0,1,\dots$, the inequality 
\begin{equation}\label{202012171848}
N(u_1^{(m)},u_2^{(m)})
\leq
C\varepsilon
+
C
N(u_1^{(m-1)},u_2^{(m-1)})^2
+
CN(u_1^{(m-1)},u_2^{(m-1)})^3
\end{equation}
holds.
\end{proposition}
{\it Proof of Proposition $\ref{20201030}$}. 
The next two lemmas will be used in the proof. 
The first lemma states that 
the null condition is preserved 
under the differentiation. 
\begin{lemma}
Suppose that 
$\{A^{\alpha\beta}:\alpha,\beta=0,1,2,3\}$ satisfies 
the null condition, that is to say, 
$A^{\alpha\beta}X_\alpha X_\beta\equiv 0$ 
for $X\in{\mathbb R}^4$ satisfying $X_0^2=X_1^2+X_2^2+X_3^2$. 
Then, for any $\Gamma_l$ $(l=0,\dots,10)$ we have
\begin{align}
\Gamma_l&
A^{\alpha\beta}
(\partial_\alpha v)
(\partial_\beta w)\\
&
=
A^{\alpha\beta}
(\partial_\alpha \Gamma_l v)
(\partial_\beta w)
+
A^{\alpha\beta}
(\partial_\alpha v)
(\partial_\beta \Gamma_l w)
+
{\tilde A}_l^{\alpha\beta}
(\partial_\alpha v)
(\partial_\beta w),\nonumber
\end{align}
with the set of the new coefficients 
$\{{\tilde A}_l^{\alpha\beta}\}$ also satisfying the null condition. 
\end{lemma} 
See \cite[Lemma 6.6.5]{Hor} for the proof. 

The second lemma says that 
the null condition creates cancellation which 
allows us to handle the quadratic nonlinear term 
as higher-order one in terms of time decay 
if its coefficients satisfy the null condition. 
\begin{lemma}\label{202012281819}
If $\{A^{\alpha\beta}\}$ satisfies the null condition, 
then 
\begin{equation}\label{202012211743}
|
A^{\alpha\beta}
(\partial_\alpha v)
(\partial_\beta w)
|
\leq
C(1+t)^{-1}
\sum_{|a|=1}
\bigl(
|\Gamma^a v|
|\partial w|
+
|\partial v|
|\Gamma^a w|
\bigr)
\end{equation}
holds.  
\end{lemma}
See, e.g., \cite[pp.\,90--91]{Al2010}, for the proof. 

Let us start with 
the estimate of $W(u_1^{(m)})$ and $W(u_2^{(m)})$. 
It suffices to discuss how to estimate the latter, 
because we can handle the former in the same way. 
\subsection{}Let us begin with the estimate of 
$\||D|^{1/6}\Gamma^\alpha u_2^{(m)}(t)\|_2$, 
$|\alpha|\leq 2$. 
Using the Li-Zhou inequality (\ref{202012211613}), we get
\begin{align}
\|
&|D|^{1/6}\Gamma^{\alpha}u_2^{(m)}(t)
\|_2\\
&
\leq
\|
|D|^{1/6}
\bigl((\Gamma^{\alpha}u_2^{(m)})(0)\bigr)
\|_2
+
\|
|D|^{-5/6}
\bigl((\partial_t\Gamma^{\alpha}u_2^{(m)})(0)\bigr)
\|_2
\nonumber\\
&
+
C
\int_0^t
\|
A_2^{11,\alpha\beta}
(\partial_\alpha u_1^{(m-1)}(\tau))
(\partial_\beta u_1^{(m-1)}(\tau))
\|_{\Gamma,2,p_1,\chi_1}d\tau\nonumber\\
&
+
C
\int_0^t
\|
A_2^{12,\alpha\beta}
(\partial_\alpha u_1^{(m-1)}(\tau))
(\partial_\beta u_2^{(m-1)}(\tau))
\|_{\Gamma,2,p_1,\chi_1}d\tau\nonumber\\
&
+
C
\int_0^t
\|
A_2^{22,\alpha\beta}
(\partial_\alpha u_2^{(m-1)}(\tau))
(\partial_\beta u_2^{(m-1)}(\tau))
\|_{\Gamma,2,p_1,\chi_1}d\tau\nonumber\\
&
+C\int_0^t
\langle\tau\rangle^{-2/3}
\|
A_2^{11,\alpha\beta}
(\partial_\alpha u_1^{(m-1)}(\tau))
(\partial_\beta u_1^{(m-1)}(\tau))
\|_{\Gamma,2,1,p_2,\chi_2}d\tau\nonumber\\
&
+C\int_0^t
\langle\tau\rangle^{-2/3}
\|
A_2^{12,\alpha\beta}
(\partial_\alpha u_1^{(m-1)}(\tau))
(\partial_\beta u_2^{(m-1)}(\tau))
\|_{\Gamma,2,1,p_2,\chi_2}d\tau\nonumber\\
&
+C\int_0^t
\langle\tau\rangle^{-2/3}
\|
A_2^{22,\alpha\beta}
(\partial_\alpha u_2^{(m-1)}(\tau))
(\partial_\beta u_2^{(m-1)}(\tau))
\|_{\Gamma,2,1,p_2,\chi_2}d\tau\nonumber\\
&
+
C
\int_0^t
\|
C_2(u^{(m-1)}(\tau),\partial u^{(m-1)}(\tau))
\|_{\Gamma,2,p_1,\chi_1}d\tau\nonumber\\
&
+
C\int_0^t
\langle\tau\rangle^{-2/3}
\|
C_2(u^{(m-1)}(\tau),\partial u^{(m-1)}(\tau))
\|_{\Gamma,2,1,p_2,\chi_2}d\tau.\nonumber
\end{align}
Here $1/p_1=1/2+5/18$, $1/p_2=1/2+1/6$, 
and we have set $A_2^{21,\alpha\beta}=0$ 
$(\alpha,\beta=0,1,2,3)$ without loss of generality. 
In the following discussion, we sometimes drop the superscript 
$(m-1)$ for simplicity. 
Using the Klainerman-Sobolev inequality (\ref{202012211620}), 
we get
\begin{align}\label{202011251636}
\|&
A_2^{11,\alpha\beta}
(\partial_\alpha u_1)
(\partial_\beta u_1)
\|_{\Gamma,2,p_1,\chi_1}
\leq
C
\|
\partial u_1
\|_{\Gamma,1,18/5,\chi_1}
\|
\partial u_1
\|_{\Gamma,2,2}\\
&
\leq
C
\langle\tau\rangle^{-3(1/2-5/18)}
\|
\partial u_1
\|_{\Gamma,2,2}^2
\leq
C\langle\tau\rangle^{-2/3}X(u_1)^2.\nonumber
\end{align}
Moreover, owing to the well-known inequality
\begin{equation}
(1+|t-r|)
|\partial v(t,x)|
\leq
C
\sum_{|\alpha|\leq 1}
|\Gamma^\alpha v(t,x)|
\end{equation}
which is a direct consequence of (\ref{20210124926}) 
and its analogue for $\partial_t$ 
(see \cite[p.\,115]{Kl87}), together with 
the Klainerman-Sobolev inequality (\ref{202012211620}), 
we get
\begin{align}\label{202010301716}
\|&
A_2^{12,\alpha\beta}
(\partial_\alpha u_1)
(\partial_\beta u_2)
\|_{\Gamma,2,p_1,\chi_1}\\
&
\leq
C\langle\tau\rangle^{-1}
\sum_{(i,j)=(1,2),(2,1)}
\|
u_i
\|_{\Gamma,3,2}
\|
\partial u_j
\|_{\Gamma,2,18/5,\chi_1}\nonumber\\
&
\leq
C\langle\tau\rangle^{-1+(1/2+2\delta)-(2/3-\delta)}
\bigl(
W(u_1)+W(u_2)
\bigr)
\bigl(
X(u_1)+Y(u_2)
\bigr)\nonumber\\
&
=
C\langle\tau\rangle^{-7/6+3\delta}
\bigl(
W(u_1)+W(u_2)
\bigr)
\bigl(
X(u_1)+Y(u_2)
\bigr).\nonumber
\end{align}
Similarly, we obtain 
\begin{equation}\label{202010301717}
\|
A_2^{22,\alpha\beta}
(\partial_\alpha u_2)
(\partial_\beta u_2)
\|_{\Gamma,2,p_1,\chi_1}
\leq
C\langle\tau\rangle^{-7/6+3\delta}
W(u_2)Y(u_2).
\end{equation}
For the estimate of the norm over $\{x:|x|>(1+\tau)/2\}$, 
we get
\begin{align}\label{202011251639}
\langle&\tau\rangle^{-2/3}
\|
A_2^{11,\alpha\beta}
(\partial_\alpha u_1)
(\partial_\beta u_1)
\|_{\Gamma,2,1,p_2,\chi_2}
\leq
C
\langle\tau\rangle^{-2/3}
\|\partial u_1\|_{\Gamma,1,2,6}
\|\partial u_1\|_{\Gamma,2,2}\\
&
\leq
C\langle\tau\rangle^{-2/3}
\|\partial u_1\|_{\Gamma,2,2}^2
\leq
C\langle\tau\rangle^{-2/3}
X(u_1)^2,\nonumber
\end{align}
where we have employed the Sobolev embedding on $S^2$. 
Moreover, by virtue of the null condition (\ref{Aass2}) 
and the Sobolev embedding on $S^2$, 
we get
\begin{align}\label{202011181806}
\langle&\tau\rangle^{-2/3}
\|
A_2^{12,\alpha\beta}
(\partial_\alpha u_1)
(\partial_\beta u_2)
\|_{\Gamma,2,1,p_2,\chi_2}\\
&
\leq
C
\langle\tau\rangle^{-5/3}
\sum_{(i,j)=(1,2),(2,1)}
\|u_i\|_{\Gamma,3,2}
\|\partial u_j\|_{\Gamma,2,2,6}\nonumber\\
&
\leq
C
\langle\tau\rangle^{-5/3}
\sum_{(i,j)=(1,2),(2,1)}
\|u_i\|_{\Gamma,3,2}
\|\partial u_j\|_{\Gamma,3,2}\nonumber\\
&
\leq
C
\langle\tau\rangle^{-7/6+3\delta}
\bigl(W(u_1)+W(u_2)\bigr)
\bigl(X(u_1)+Y(u_2)\bigr).\nonumber
\end{align}
Here, we have used (\ref{202012211743}). 
Similarly, we get thanks to the null condition (\ref{Aass3}) 
\begin{equation}\label{20210124940}
\langle\tau\rangle^{-2/3}
\|
A_2^{22,\alpha\beta}
(\partial_\alpha u_2)
(\partial_\beta u_2)
\|_{\Gamma,2,1,p_2,\chi_2}
\leq
C
\langle\tau\rangle^{-7/6+3\delta}
W(u_2)Y(u_2).
\end{equation}
Note that in the discussion 
(\ref{202011181806})--(\ref{20210124940}), 
we have relied upon the null condition (\ref{Aass2})--(\ref{Aass3}) 
so as to beat the growth of 
$\|\partial u_2(t)\|_{\Gamma,3,2}$ 
(see the definition of $Y(v)$, (\ref{DefY(v)})). 

Let us turn our attention to the cubic terms. 
Using the inequality 
$\||D|^{1/6}\partial v\|_2\leq\|\partial v\|_{H^1}$ 
and modifying the argument 
in (\ref{202010301816}) and (\ref{202010301840}) 
suitably, we obtain for $i=1,2$
\begin{align}\label{202011261048}
\|&
C_i(u,\partial u)
\|_{\Gamma,2,p_1,\chi_1}
+
\langle\tau\rangle^{-2/3}
\|
C_i(u,\partial u)
\|_{\Gamma,2,1,p_2,\chi_2}\\
&
\leq
C\sum_{l=0}^3
\langle\tau\rangle^{-5/3+(3-l)/3+l\delta}
\bigl(
W(u_1)+W(u_2)
\bigr)^{3-l}
\bigl(
X(u_1)+Y(u_2)
\bigr)^l\nonumber\\
&
\leq
C\langle\tau\rangle^{-2/3}
\sum_{l=0}^3
\bigl(
W(u_1)+W(u_2)
\bigr)^{3-l}
\bigl(
X(u_1)+Y(u_2)
\bigr)^l.\nonumber
\end{align}
We have finished the estimate 
$\||D|^{1/6}\Gamma^\alpha u_2^{(m)}(t)\|_2$, 
$|\alpha|\leq 2$. 
\subsection{}Note that in (\ref{202010301716}) 
and (\ref{202011181806}), 
we have used the estimation lemma (\ref{202012211743}) 
and thereby gained more time decay, 
at the expense of a kind of loss of derivatives; 
namely, the norm $\|u_i\|_{\Gamma,3,2}$ has appeared. 
In order to close the chain of estimates, 
let us next discuss how to bound it. 

When estimating $\Gamma^\alpha u_2^{(m)}(t)$ 
for $|\alpha|=3$, 
we can afford no loss of derivatives. 
Therefore, we can no longer rely upon the estimation lemma 
(\ref{202012211743}). 
This is the reason why 
we have only a coarse estimate 
for $\|u_2^{(m)}(t)\|_{\Gamma,3,2}$. 
For our aim, however, it is sufficient. 
Employing the argument in (\ref{202011251636}) and 
(\ref{202011251639}), where we have never relied upon 
the null condition, and using 
$\|\cdots\|_{\Gamma,3,6/5,\chi_1}$ 
and 
$\langle\tau\rangle^{-1/2}\|\cdots\|_{\Gamma,3,1,4/3,\chi_2}$ 
in place of 
$\|\cdots\|_{\Gamma,2,p_1,\chi_1}$ 
and 
$\langle\tau\rangle^{-2/3}\|\cdots\|_{\Gamma,2,1,p_2,\chi_2}$, 
respectively, we get 
\begin{align}\label{202012111737}
\|&
A_2^{ij,\alpha\beta}
(\partial_\alpha u_i)
(\partial_\beta u_j)
\|_{\Gamma,3,6/5,\chi_1}
+
\langle\tau\rangle^{-1/2}\|A_2^{ij,\alpha\beta}
(\partial_\alpha u_i)
(\partial_\beta u_j)\|_{\Gamma,3,1,4/3,\chi_2}\\
&
\leq
C\langle\tau\rangle^{-1/2+2\delta}
\bigl(
X(u_1)+Y(u_2)
\bigr)^2.\nonumber
\end{align}
Particular attention should be paid for the cubic terms, 
because we can never rely upon 
the norm 
$\||D|^{1/6}\Gamma^\alpha u(t)\|_2$ 
$(|\alpha|=3)$, and therefore 
it is necessary to handle them 
differently from (\ref{202011261048}). 
Using the H\"older inequality and 
the Klainerman-Sobolev inequality (\ref{202012211620}), 
and the Sobolev embedding ${\dot H}^{1/6}\hookrightarrow L^{9/4}$, 
we get 
\begin{align}\label{202011261441}
\|&u_iu_ju_k\|_{\Gamma,3,6/5,\chi_1}\\
&
\leq
\|u_i\|_{\Gamma,1,6,\chi_1}
\|u_j\|_{\Gamma,1,6,\chi_1}
\|u_k\|_{\Gamma,3,2}\nonumber\\
&
\leq
C\langle\tau\rangle^{-5/3+2(1/3)+(1/2+2\delta)}
\biggl(
\langle\tau\rangle^{-1/3}
\sum_{|\alpha|\leq 2}
\||D|^{1/6}\Gamma^\alpha u_i\|_2
\biggr)\nonumber\\
&
\times
\biggl(
\langle\tau\rangle^{-1/3}
\sum_{|\alpha|\leq 2}
\||D|^{1/6}\Gamma^\alpha u_j\|_2
\biggr)
\biggl(
\langle\tau\rangle^{-1/2-2\delta}
\|u_k\|_{\Gamma,3,2}
\biggr)\nonumber\\
&
\leq
C\langle\tau\rangle^{-1/2+2\delta}
\bigl(
W(u_1)+W(u_2)
\bigr)^3.\nonumber
\end{align}
Moreover, using the inequality 
$\||D|^{1/6}\partial v\|_{L^2}\leq\|\partial v\|_{H^1}$ 
and repeating the same argument as in (\ref{202011261441}), 
we get
\begin{align}\label{202011261500}
\|&u_iu_j\partial u_k\|_{\Gamma,3,6/5,\chi_1}
+
\|u_i(\partial u_j)(\partial u_k)\|_{\Gamma,3,6/5,\chi_1}
+
\|(\partial u_i)(\partial u_j)
(\partial u_k)\|_{\Gamma,3,6/5,\chi_1}\\
&
\leq
C\langle\tau\rangle^{-5/6+3\delta}
\sum_{l=1}^3
\bigl(
W(u_1)+W(u_2)
\bigr)^{3-l}
\bigl(
X(u_1)+Y(u_2)
\bigr)^l.\nonumber
\end{align}
Combining (\ref{202011261441})--(\ref{202011261500}), 
we have obtained 
\begin{equation}
\|C_i(u,\partial u)\|_{\Gamma,3,6/5,\chi_1}
\leq
C\langle\tau\rangle^{-1/2+2\delta}
\sum_{l=0}^3
\bigl(
W(u_1)+W(u_2)
\bigr)^{3-l}
\bigl(
X(u_1)+Y(u_2)
\bigr)^l.
\end{equation}
Furthermore, 
using the H\"older inequality, 
the Sobolev embedding on $S^2$, 
and the weighted inequality (\ref{202012211106}), 
we get
\begin{align}\label{202011261552}
\langle&\tau\rangle^{-1/2}
\|u_iu_ju_k\|_{\Gamma,3,1,4/3,\chi_2}\\
&
\leq
\langle\tau\rangle^{-1/2}
\|u_i\|_{\Gamma,1,4,8,\chi_2}
\|u_j\|_{\Gamma,1,4,8,\chi_2}
\|u_k\|_{\Gamma,3,2}\nonumber\\
&
\leq
C\langle\tau\rangle^{-1/2}
\|u_i\|_{\Gamma,2,4,2,\chi_2}
\|u_j\|_{\Gamma,2,4,2,\chi_2}
\|u_k\|_{\Gamma,3,2}\nonumber\\
&
\leq
C\langle\tau\rangle^{-1/2-2(1/2-1/4)\times 2}
\biggl(
\sum_{\substack{|\alpha|\leq 2\\i=1,2}}
\||D|^{1/4}\Gamma^\alpha u_i\|_2
\biggr)^2
\|u_k\|_{\Gamma,3,2}\nonumber\\
&
\leq
C\langle\tau\rangle^{-3/2+2(1/4)+(1/2+2\delta)}
\biggl(
\sum_{\substack{|\alpha|\leq 2\\i=1,2}}
\langle\tau\rangle^{-1/4}
\||D|^{1/4}\Gamma^\alpha u_i\|_2
\biggr)^2
\bigl(
\langle\tau\rangle^{-1/2-2\delta}
\|u_k\|_{\Gamma,3,2}
\bigr)\nonumber\\
&
\leq
C\langle\tau\rangle^{-1/2+2\delta}
\bigl(
W(u_1)+W(u_2)
\bigr)^3.\nonumber
\end{align}
We must notice that 
the norm $\||D|^{1/4}\Gamma^\alpha u_i\|_2$ 
$(|\alpha|\leq 2)$, which will be estimated 
in the next subsection, has appeared here. 

Moreover, using the inequality 
$\||D|^{1/4}\partial v\|_{L^2}\leq\|\partial v\|_{H^1}$ 
and repeating the same argument as in (\ref{202011261552}), 
we get
\begin{align}\label{202011261645}
\langle&\tau\rangle^{-1/2}
\bigl(\|u_iu_j\partial u_k\|_{\Gamma,3,1,4/3,\chi_2}\\
&
\hspace{2cm}
+
\|u_i(\partial u_j)(\partial u_k)\|_{\Gamma,3,1,4/3,\chi_2}
+
\|(\partial u_i)(\partial u_j)(\partial u_k)\|_{\Gamma,3,1,4/3,\chi_2}
\bigr)\nonumber\\
&
\leq
C\langle\tau\rangle^{-3/4+3\delta}
\sum_{l=1}^3
\bigl(
W(u_1)+W(u_2)
\bigr)^{3-l}
\bigl(
X(u_1)+Y(u_2)
\bigr)^l.\nonumber
\end{align}
Combining (\ref{202011261552}) and (\ref{202011261645}), 
we have shown
\begin{align}\label{202012111739}
\langle&\tau\rangle^{-1/2}
\|C_i(u,\partial u)\|_{\Gamma,3,1,4/3,\chi_2}\\
&
\leq
C\langle\tau\rangle^{-1/2+2\delta}
\sum_{l=0}^3
\bigl(
W(u_1)+W(u_2)
\bigr)^{3-l}
\bigl(
X(u_1)+Y(u_2)
\bigr)^l.\nonumber
\end{align}
By (\ref{202012111737})--(\ref{202012111739}), 
we have shown for $|\alpha|\leq 3$
\begin{align}
\langle&t\rangle^{-1/2-2\delta}
\|
\Gamma^\alpha u_2^{(m)}(t)
\|_2\\
&
\leq
\|
(\Gamma^\alpha u_2^{(m)})(0)
\|_2
+
\|
|D|^{-1}
\bigl(
(\partial_t\Gamma^\alpha u_2^{(m)})(0)
\bigr)
\|_2
\nonumber\\
&
+
C
\bigl(
X(u_1^{(m-1)})
+
Y(u_2^{(m-1)})
\bigr)^2\nonumber\\
&
+
C
\sum_{l=0}^3
\bigl(
W(u_1^{(m-1)})+W(u_2^{(m-1)})
\bigr)^{3-l}
\bigl(
X(u_1^{(m-1)})+Y(u_2^{(m-1)})
\bigr)^l.\nonumber
\end{align}
\subsection{}Next, we must consider the estimate of 
$\||D|^{1/4}\Gamma^\alpha u_2^{(m)}(t)\|_2$, 
$|\alpha|\leq 2$ which has appeared in (\ref{202011261552}) above. 
In fact, we have only to repeat essentially the same argument 
as we have done above for the estimate of 
$\||D|^{1/6}\Gamma^\alpha u_2^{(m)}(t)\|_2$, 
$|\alpha|\leq 2$. 
Indeed, using appropriately the norms $\|\cdots\|_{\Gamma,2,4/3,\chi_1}$ and 
$\langle\tau\rangle^{-3/4}\|\cdots\|_{\Gamma,2,1,8/5,\chi_2}$ 
in place of $\|\cdots\|_{\Gamma,2,p_1,\chi_1}$ 
and $\langle\tau\rangle^{-2/3}\|\cdots\|_{\Gamma,2,1,p_2,\chi_2}$ 
respectively, 
we obtain for $|\alpha|\leq 2$
\begin{align}
\langle&t\rangle^{-1/4}
\|
|D|^{1/4}
\Gamma^\alpha u_2^{(m)}(t)
\|_2\\
&
\leq
\|
|D|^{1/4}
\bigl(
(\Gamma^\alpha u_2^{(m)})(0)
\bigr)
\|_2
+
\|
|D|^{-3/4}
\bigl(
(\partial_t\Gamma^\alpha u_2^{(m)})(0)
\bigr)
\|_2
\nonumber\\
&
+
CX(u_1^{(m-1)})^2
+
C
\bigl(
W(u_1^{(m-1)})+W(u_2^{(m-1)})
\bigr)
\bigl(
X(u_1^{(m-1)})+Y(u_2^{(m-1)})
\bigr)\nonumber\\
&
+
C\sum_{l=0}^3
\bigl(
W(u_1^{(m-1)})+W(u_2^{(m-1)})
\bigr)^{3-l}
\bigl(
X(u_1^{(m-1)})+Y(u_2^{(m-1)})
\bigr)^l.\nonumber
\end{align}
\subsection{}Next, our concern centers on 
the estimate of $X(u_1^{(m)})$ and $Y(u_2^{(m)})$. 
Dealing with the quadratic nonlinear terms as in 
\cite{Al2010} and \cite{HY2017}, 
we obtain
\begin{align}\label{202012231839}
&
\sum_{|\alpha|\leq 3}
\bigl(
\|
\partial\Gamma^\alpha u_1^{(m)}(t)
\|_2
+
\|
\langle
\tau-|x|
\rangle^{-1/2-\eta}
T\Gamma^\alpha u_1^{(m)}
\|_{L^2((0,t)\times{\mathbb R}^3)}
\bigr)\\
&
+
\sum_{|\alpha|\leq 3}
\bigl(
\langle t\rangle^{-\delta}
\|
\partial\Gamma^\alpha u_2^{(m)}(t)
\|_2
+
\langle t\rangle^{-\delta}
\|
\langle
\tau-|x|
\rangle^{-1/2-\eta}
T\Gamma^\alpha u_2^{(m)}
\|_{L^2((0,t)\times{\mathbb R}^3)}
\bigr)\nonumber\\
&
\leq
C\sum_{\substack{i=1,2\\|\alpha|\leq 3}}
\|
(\partial\Gamma^\alpha u_i^{(m)})(0)
\|_2
+
C
\bigl(
X(u_1^{(m-1)})
+
Y(u_2^{(m-1)})
\bigr)^2\nonumber\\
&
+
C\int_0^t
\|
C_1(u^{(m-1)},\partial u^{(m-1)})
\|_{\Gamma,3,2}
d\tau
+
C
\int_0^t
\|
C_2(u^{(m-1)},\partial u^{(m-1)})
\|_{\Gamma,3,2}
d\tau.\nonumber
\end{align}
For readers' convenience, 
we give the outline of the proof of (\ref{202012231839}) 
in Appendix. 
As above, we must pay careful enough attention to 
the cubic terms, 
and we handle 
$\|C_i(u^{(m-1)},\partial u^{(m-1)})\|_{\Gamma,3,2,\chi_1}$ 
and $\|C_i(u^{(m-1)},\partial u^{(m-1)})\|_{\Gamma,3,2,\chi_2}$, 
separately. 

Again, we will drop the superscript 
$(m-1)$ for simplicity. 
Using the Klainerman-Sobolev inequality (\ref{202012211648}), 
we get
\begin{align}
\sum_{i,j,k=1}^2&
\|u_iu_ju_k\|_{\Gamma,3,2,\chi_1}
\leq
C
\sum_{i,j,k=1}^2
\|u_i\|_{\Gamma,1,\infty,\chi_1}
\|u_j\|_{\Gamma,1,\infty,\chi_1}
\|u_k\|_{\Gamma,3,2}\\
&
\leq
C
\langle\tau\rangle^{-3}
\sum_{i,j,k=1}^2
\|u_i\|_{\Gamma,3,2}
\|u_j\|_{\Gamma,3,2}
\|u_k\|_{\Gamma,3,2}\nonumber\\
&
\leq
C\langle\tau\rangle^{-3/2+6\delta}
\bigl(
W(u_1)+W(u_2)
\bigr)^3.\nonumber
\end{align}
We treat the other cubic terms similarly, to get
\begin{align}
\|&
C_i(u,\partial u)
\|_{\Gamma,3,2,\chi_1}\\
&
\leq
C\sum_{l=0}^3
\langle\tau\rangle^{-3/2+6\delta-l/2-l\delta}
\bigl(
W(u_1)+W(u_2)
\bigr)^{3-l}
\bigl(
X(u_1)+Y(u_2)
\bigr)^l.\nonumber
\end{align}
As for the estimate of $\|C_i(u,\partial u)\|_{\Gamma,3,2,\chi_2}$, 
it is necessary to proceed in a completely different way. 
It is obvious that we need to handle the following 4 types of norms: 
(i) $\|u_iu_ju_k\|_{\Gamma,3,2,\chi_2}$, 
(ii) $\|u_iu_j\partial u_k\|_{\Gamma,3,2,\chi_2}$, 
(iii) $\|u_i(\partial u_j)(\partial u_k)\|_{\Gamma,3,2,\chi_2}$, 
and 
(iv) $\|(\partial u_i)(\partial u_j)(\partial u_k)\|_{\Gamma,3,2,\chi_2}$, 
which will be estimated separately. 

\noindent $\cdot$ Bound for (i). 
Using the Sobolev embedding on $S^2$ 
and the weighted inequality (\ref{202012211106}), we proceed as 
\begin{align}
&
\sum_{|\alpha|,|\beta|,|\gamma|\leq 1}
\|(\Gamma^\alpha u_i)(\Gamma^\beta u_j)
(\Gamma^\gamma u_k)\|_{2,\chi_2}\\
&
+
\sum_{|\alpha|\leq 2,|\beta|\leq 1}
\|(\Gamma^\alpha u_i)(\Gamma^\beta u_j)u_k\|_{2,\chi_2}
+
\sum_{|\alpha|\leq 3}
\|(\Gamma^\alpha u_i)u_ju_k\|_{2,\chi_2}\nonumber\\
&
\leq
\sum_{|\alpha|,|\beta|,|\gamma|\leq 1}
\|\Gamma^\alpha u_i\|_{6,\chi_2}
\|\Gamma^\beta u_j\|_{6,\chi_2}
\|\Gamma^\gamma u_k\|_{6,\chi_2}\nonumber\\
&
+
\sum_{|\alpha|\leq 2,|\beta|\leq 1}
\|\Gamma^\alpha u_i\|_{6,4,\chi_2}
\|\Gamma^\beta u_j\|_{6,4,\chi_2}
\|u_k\|_{6,\infty,\chi_2}\nonumber\\
&
+
\sum_{|\alpha|\leq 3}
\|\Gamma^\alpha u_i\|_{6,2,\chi_2}
\|u_j\|_{6,\infty,\chi_2}
\|u_k\|_{6,\infty,\chi_2}\nonumber\\
&
\leq
\sum_{|\alpha|\leq 3,\,|\beta|,|\gamma|\leq 2}
\|\Gamma^\alpha u_i\|_{6,2,\chi_2}
\|\Gamma^\beta u_j\|_{6,2,\chi_2}
\|\Gamma^\gamma u_k\|_{6,2,\chi_2}\nonumber\\
&
\leq
C\langle\tau\rangle^{-2}
\sum_{|\alpha|\leq 3,\,|\beta|,|\gamma|\leq 2}
\||D|^{1/3}\Gamma^\alpha u_i\|_2
\||D|^{1/3}\Gamma^\beta u_j\|_2
\||D|^{1/3}\Gamma^\gamma u_k\|_2.\nonumber
\end{align}
Using for $|\alpha|\leq 3$ and $|\beta|\leq 2$
\begin{align}
&
\|
|D|^{1/3}\Gamma^\alpha u_i
\|_2
\leq
\|
\Gamma^\alpha u_i
\|_2^{2/3}
\|
|D|\Gamma^\alpha u_i
\|_2^{1/3},\label{2021interpolation1}\\
&
\|
|D|^{1/3}\Gamma^\beta u_j
\|_2
\leq
\|
|D|^{1/6}
\Gamma^\beta u_j
\|_2^{4/5}
\|
|D|\Gamma^\beta u_j
\|_2^{1/5},\label{2021interpolation2}
\end{align}
we obtain
\begin{align}
\|&u_iu_ju_k\|_{\Gamma,3,2,\chi_2}\\
&
\leq
C
\langle\tau\rangle^{-17/15+31\delta/15}
\bigl(
W(u_1)+W(u_2)
\bigr)^{34/15}
\bigl(
X(u_1)+Y(u_2)
\bigr)^{11/15}.\nonumber
\end{align}

\noindent$\cdot$ Bound for (ii). We deal with 
$\|\bigl(\Gamma^\alpha (u_iu_j)\bigr)\partial u_k\|_{2,\chi_2}$ 
$(|\alpha|\leq 3)$ 
and 
$\|\bigl(\Gamma^\alpha (u_iu_j)\bigr)
\Gamma^\beta\partial u_k\|_{2,\chi_2}$ 
$(|\alpha|+|\beta|\leq 3,\,|\alpha|\leq 2)$, 
separately. For the former, we proceed 
by using the Sobolev embedding on $S^2$, 
(\ref{20201221938}), and (\ref{20201221933})
\begin{align}
&
\sum_{i,j,k=1}^2
\sum_{|\alpha|\leq 3}
\|\bigl(\Gamma^\alpha (u_iu_j)\bigr)\partial u_k\|_{2,\chi_2}\\
&
\leq
C
\sum_{i,j,k=1}^2
\|u_i\|_{\Gamma,1,\infty,\chi_2}
\|u_j\|_{\Gamma,3,2}
\|\partial u_k\|_{\infty,\chi_2}\nonumber\\
&
\leq
C\sum_{i,j,k=1}^2
\|u_i\|_{\Gamma,2,\infty,24/11,\chi_2}
\|u_j\|_{\Gamma,3,2}
\|\partial u_k\|_{\Gamma,1,\infty,4,\chi_2}\nonumber\\
&
\leq
C
\langle\tau\rangle^{-3/2+7/12}\langle\tau\rangle^{-1}
\sum_{i,j,k=1}^2
\biggl(
\sum_{|\alpha|\leq 2}
\||D|^{7/12}\Gamma^\alpha u_i\|_2
\biggr)
\|u_j\|_{\Gamma,3,2}
\|\partial u_k\|_{\Gamma,2,2}\nonumber\\
&
\leq
C
\langle\tau\rangle^{-23/12+(1/2+2\delta)+\delta}
\biggl(
\sum_{i=1}^2
\sum_{|\alpha|\leq 2}
\||D|^{1/6}\Gamma^\alpha u_i\|_2^{1/2}
\||D|\Gamma^\alpha u_i\|_2^{1/2}
\biggr)
\nonumber\\
&
\times
\bigl(
W(u_1)+W(u_2)
\bigr)
\bigl(
X(u_1)+X(u_2)
\bigr)\nonumber\\
&
\leq
C
\langle\tau\rangle^{-5/4+7\delta/2}
\bigl(
W(u_1)+W(u_2)
\bigr)^{3/2}
\bigl(
X(u_1)+X(u_2)
\bigr)^{3/2}.\nonumber
\end{align}
For the latter, we get by the Sobolev embedding on $S^2$ 
and (\ref{20201221938}), 
\begin{align}
&
\sum_{i,j,k=1}^2
\sum_{\substack{|\alpha|+|\beta|\leq 3\\|\alpha|\leq 2}}
\|\bigl(\Gamma^\alpha (u_iu_j)\bigr)
\Gamma^\beta\partial u_k\|_{2,\chi_2}\\
&
\leq
C
\sum_{i,j,k=1}^2
\|u_i\|_{\Gamma,1,\infty,\chi_2}
\|u_j\|_{\Gamma,2,\infty,\chi_2}
\|\partial u_k\|_{\Gamma,3,2}\nonumber\\
&
\leq
C\sum_{i,j,k=1}^2
\|u_i\|_{\Gamma,2,\infty,24/11,\chi_2}
\|u_j\|_{\Gamma,3,\infty,4,\chi_2}
\|\partial u_k\|_{\Gamma,3,2}\nonumber\\
&
\leq
C\langle\tau\rangle^{-3/2+7/12}\langle\tau\rangle^{-1/2}\nonumber\\
&
\times
\sum_{i,j,k=1}^2
\biggl(
\sum_{|\alpha|\leq 2}
\||D|^{7/12}\Gamma^\alpha u_i\|_2
\biggr)
\biggl(
\sum_{|\alpha|\leq 3}
\||D|\Gamma^\alpha u_j\|_2
\biggr)
\|\partial u_k\|_{\Gamma,3,2}\nonumber\\
&
\leq
C\langle\tau\rangle^{-3/2+7/12}\langle\tau\rangle^{-1/2}
\langle\tau\rangle^{2\delta}\nonumber\\
&
\times
\biggl(
\sum_{i=1}^2
\sum_{|\alpha|\leq 2}
\||D|^{1/6}\Gamma^\alpha u_i\|_2^{1/2}
\||D|\Gamma^\alpha u_i\|_2^{1/2}
\biggr)
\bigl(
X(u_1)+Y(u_2)
\bigr)^2\nonumber\\
&
\leq
C\langle\tau\rangle^{-5/4+5\delta/2}
\bigl(
W(u_1)+W(u_2)
\bigr)^{1/2}
\bigl(
X(u_1)+Y(u_2)
\bigr)^{5/2}.\nonumber
\end{align}

\noindent$\cdot$ Bound for (iii). 
We handle the norms of this type 
by estimating 
$\|(\Gamma^\alpha u_i)(\partial u_j)(\partial u_k)\|_{2,\chi_2}$ 
$(|\alpha|\leq 3)$
and the others 
$\|(\Gamma^\alpha u_i)
\Gamma^\beta
\bigl(
(\partial u_j)(\partial u_k)
\bigr)\|_{2,\chi_2}$ 
$(|\alpha|+|\beta|\leq 3,\,|\alpha|\leq 2)$, 
separately. 
The former is estimated as 
\begin{align}
&\sum_{i,j,k=1}^2
\sum_{|\alpha|\leq 3}
\|(\Gamma^\alpha u_i)(\partial u_j)(\partial u_k)\|_{2,\chi_2}\\
&
\leq
\sum_{i,j,k=1}^2
\sum_{|\alpha|\leq 3}
\|u_i\|_{\Gamma,3,2}
\|\partial u_j\|_{\infty,\chi_2}
\|\partial u_k\|_{\infty,\chi_2}\nonumber\\
&
\leq
C\langle\tau\rangle^{1/2+2\delta}
\bigl(
W(u_1)+W(u_2)
\bigr)
\sum_{j,k=1}^2
\|\partial u_j\|_{\Gamma,1,\infty,4,\chi_2}
\|\partial u_k\|_{\Gamma,1,\infty,4,\chi_2}\nonumber\\
&
\leq
C\langle\tau\rangle^{-3/2+4\delta}
\bigl(W(u_1)+W(u_2)\bigr)
\bigl(X(u_1)+Y(u_2)\bigr)^2,\nonumber
\end{align}
where we have used the Sobolev embedding on $S^2$ 
and (\ref{20201221933}). Moreover, we handle the latter as 
\begin{align}
&\sum_{i,j,k=1}^2
\sum_{\substack{|\alpha|+|\beta|\leq 3\\|\alpha|\leq 2}}
\|(\Gamma^\alpha u_i)
\Gamma^\beta
\bigl(
(\partial u_j)(\partial u_k)
\bigr)\|_{2,\chi_2}\\
&
\leq
C
\sum_{i,j,k=1}^2
\|u_i\|_{\Gamma,2,\infty,\chi_2}
\|\partial u_j\|_{\Gamma,1,\infty,\chi_2}
\|\partial u_k\|_{\Gamma,3,2}\nonumber\\
&
\leq
C
\langle\tau\rangle^\delta
\bigl(
X(u_1)+Y(u_2)
\bigr)
\sum_{i,j=1}^2
\|u_i\|_{\Gamma,3,\infty,4,\chi_2}
\|\partial u_j\|_{\Gamma,2,\infty,4,\chi_2}
\nonumber\\
&
\leq
C
\langle\tau\rangle^{-3/2+3\delta}
\bigl(
X(u_1)+Y(u_2)
\bigr)^3,\nonumber
\end{align}
where we have used the trace-type inequalities 
(\ref{20201221932})--(\ref{20201221933}).

\noindent$\cdot$ Bound for (iv). 
Using the trace-type inequality (\ref{20201221933}) suitably, 
we easily obtain
\begin{equation}
\sum_{i,j,k=1}^2
\|
(\partial u_i)
(\partial u_j)
(\partial u_k)
\|_{\Gamma,3,2,\chi_2}
\leq
C\langle\tau\rangle^{-2+3\delta}
\bigl(
X(u_1)+Y(u_2)
\bigr)^3.
\end{equation}
The estimate of $X(u_1^{(m)})$ and $Y(u_2^{(m)})$ has been finished.
\subsection{}Now, we are in a position to complete the proof of 
Proposition \ref{20201030} 
and Theorem \ref{GEforLRequation}. 
As in (\ref{202012161715}), 
it is possible to get
\begin{align}
\sum_{i=1}^2
\biggl(&
\sum_{|\alpha|\leq 2}
\bigl(
\|
|D|^{1/6}
\bigl(
(\Gamma^\alpha u_i^{(m)})(0)
\bigr)
\|_2
+
\|
|D|^{-5/6}
\bigl(
(\partial_t\Gamma^\alpha u_i^{(m)})(0)
\bigr)
\|_2
\bigr)\\
&
+
\sum_{|\alpha|\leq 3}
\bigl(
\|
(\Gamma^\alpha u_i^{(m)})(0)
\|_2
+
\|
|D|^{-1}
\bigl(
(\partial_t\Gamma^\alpha u_i^{(m)})(0)
\bigr)
\|_2
\bigr)\nonumber\\
&
+
\sum_{|\alpha|\leq 2}
\bigl(
\|
|D|^{1/4}
\bigl(
(\Gamma^\alpha u_i^{(m)})(0)
\bigr)
\|_2
+
\|
|D|^{-3/4}
\bigl(
(\partial_t\Gamma^\alpha u_i^{(m)})(0)
\bigr)
\|_2
\bigr)
\nonumber\\
&
+
\sum_{|\alpha|\leq 3}
\|
(\partial\Gamma^\alpha u_i^{(m)})(0)
\|_2
\biggr)
\leq
C\varepsilon,\nonumber
\end{align}
with a constant $C$ independent of $m$, 
provided that $\varepsilon$ is small enough. 
As remarked in the beginning of the proof of 
this proposition, 
it is possible to bound 
$W(u_1^{(m)})$ in the same way as we have done above 
for $W(u_2^{(m)})$. 
Recalling the definition of $N(u,v)$, 
we have shown the estimate (\ref{202012171848}). 

It follows from Proposition \ref{20201030} that 
there exist constants 
$0<\varepsilon_2<1$ and $C_2>0$ 
such that 
if $0<\varepsilon<\varepsilon_2$, then 
we have 
\begin{equation} 
N(u_1^{(m)},u_2^{(m)})
\leq
2C_2\varepsilon_2,
\quad
m=0,1,\dots
\end{equation}
Moreover, 
repeating essentially the same argument as 
in the proof of Proposition \ref{20201030} 
and choosing $\varepsilon_2$ smaller if necessary, 
we get for $m=2,3,\dots$,
\begin{equation}\label{202012181145}
N(u_1^{(m+1)}-u_1^{(m)},u_2^{(m+1)}-u_2^{(m)})
\leq
\biggl(
\frac12
\biggr)^{m-1}
N(u_1^{(2)}-u_1^{(1)},u_2^{(2)}-u_2^{(1)}).
\end{equation}
Equipped with the metric 
$d((u,v),(u',v'))
:=
N(u-u',v-v')$, 
the set 
${\mathcal V}(2C_2\varepsilon_2):=
\{(u,v)\in{\mathcal V}:N(u,v)\leq 2C_2\varepsilon_2\}$ 
is a complete metric space. 
(See (\ref{202012181141}) for the definition of ${\mathcal V}$.) 
Since (\ref{202012181145}) implies that 
the sequence $\{(u_1^{(m)},u_2^{(m)})\}$ defined above 
is a Cauchy sequence, 
it has a limit $(u_1,u_2)$ in ${\mathcal V}(2C_2\varepsilon_2)$. 
This is a unique solution to (\ref{ModelWeakNull}) 
that we have sought for. 
We have finished the proof of Theorem \ref{GEforLRequation}. 
\appendix
\section{}
We consider existence of global solutions to 
the Cauchy problem 
\begin{equation}\label{202101201845}
\begin{cases}
\displaystyle{
\partial_t^2 u-\Delta u
=
|v|^p,
             }
&
\displaystyle{t>0,\,x\in{\mathbb R}^n},\\
\displaystyle{
\partial_t^2 v-\Delta v
=
|\partial_t u|^{(n+1)/(n-1)},
             }&
\displaystyle{t>0,\,x\in{\mathbb R}^n},\\
\displaystyle{
u(0)=\varepsilon f_1,\,\,\partial_t u(0)=\varepsilon g_1,}
&{}\\
\displaystyle{
v(0)=\varepsilon f_2,\,\,\partial_t v(0)=\varepsilon g_2.}
&{}
\end{cases}
\end{equation}
Here, $f_i$, $g_i\in{\mathcal S}({\mathbb R}^n)$ 
$(i=1,2)$ are real-valued. 
Recall the notation $\sigma(a)=1/2-1/a$ $(a>0)$. 
We prove:
\begin{theorem}\label{202101221717}
Let $n=2,3$. 
Suppose that $p>1+3/(n-1)$. 
Then, there exist 
positive constants $C_0$ and $\varepsilon_0$ 
such that 
if $0<\varepsilon\leq\varepsilon_0$, 
then 
the Cauchy problem {\rm $($\ref{202101201845}$)$} 
admits a unique global solution $(u,v)$ satisfying 
\begin{align}
&
\esssup\displaylimits_{t>0}
(1+t)^{-1}
\|u(t)\|_2\\
&
+
\sum_{|\alpha|\leq 2}
\biggl(
\esssup\displaylimits_{t>0}
\|
\partial\Gamma^\alpha u(t)
\|_2
+
\esssup\displaylimits_{t>0}
(1+t)^{-1/(2p)}
\|
|D|^{\sigma(2p)}
\Gamma^\alpha v(t)
\|_2
\biggr)\nonumber\\
&
\leq
2C_0\varepsilon_0.\nonumber
\end{align}
\end{theorem}
\begin{remark}
By \cite[Proposition 9.1]{ISW}, 
we know that for $n=3$, 
the condition $p>5/2$ is sharp in general 
for global existence of 
small solutions to (\ref{202101201845}). 
Nonexistence of global, small solutions 
has been studied also 
in \cite{HY2016}, \cite{DFW2019} 
for systems similar to (\ref{202101201845}). 
\end{remark}

{\it Proof of Theorem $\ref{202101221717}$}. We modify the argument in 
Section \ref{section202102241600}. 
We introduce
\begin{align}
{\mathcal W}_1:=
\{
(u&,v)\in
C
(
[0,\infty)
;
H^1({\mathbb R}^n)
\times{\dot H}^{\sigma(2p)}({\mathbb R}^n)
)\,:\,\\
&
\partial_j\Gamma^\alpha u 
\in
C([0,\infty);L^2({\mathbb R}^n)),\,\,
0\leq j\leq n,\,|\alpha|\leq 1,\nonumber\\
&
|D|^{\sigma(2p)}\Gamma^\alpha v
\in
C([0,\infty);L^2({\mathbb R}^n)),\,\,|\alpha|\leq 1,\nonumber\\
&
u(0)=\varepsilon f_1,\,\partial_t u(0)=\varepsilon g_1,\,
v(0)=\varepsilon f_2,\,\partial_t v(0)=\varepsilon g_2\}\nonumber
\end{align}
and 
\begin{align}
{\mathcal W}_2:=
\{
(u&,v)\in{\mathcal W}_1\,:\,\\
&
\partial_j\Gamma^\alpha u 
\in
L^\infty
((0,\infty);L^2({\mathbb R}^n)),\,\,
0\leq j\leq n,\,|\alpha|\leq 2,\nonumber\\
&
|D|^{\sigma(2p)}\Gamma^\alpha v
\in
L^\infty
((0,\infty);L^2({\mathbb R}^n)),\,\,|\alpha|\leq 2\}.\nonumber
\end{align}
Also, for $i=1,2$ we set 
$M_i(u,v)=X_i(u)+{\hat B}_i(v)$. 
Here, 
${\hat B}_1(v)$ and ${\hat B}_2(v)$ are 
defined in the same way as in 
Section \ref{section202102241600} 
(see (\ref{hatBnorm})). 
We define $X_1(u)$ and $X_2(u)$ 
slightly differently from (\ref{normX18}). 
Namely, we set  
\begin{equation}
X_1(u)
=
\sup_{t>0}
(1+t)^{-1}
\|u(t)\|_2
+
\sum_{|\alpha|\leq 1}
\sum_{j=0}^n
\sup_{t>0}
\|
\partial_j
\Gamma^\alpha u(t)
\|_{2},
\end{equation}
and define $X_2(u)$ 
by replacing $\sup_{t>0}$ 
and $\sum_{|\alpha|\leq 1}$ 
with 
$\esssup_{t>0}$ and $\sum_{|\alpha|\leq 2}$, 
respectively. 
The set 
${\mathcal W}_1$ is complete 
with the metric 
$\eta_1
\bigl(
(u,v),(\tilde u,\tilde v)
\bigr)
:=
M_1(u-\tilde u,v-\tilde v)
$. 
Using the constants $C_0$ and $\varepsilon_0$ 
appearing below, 
we set 
${\mathcal W}_2(2C_0\varepsilon_0):=
\{(u,v)\in{\mathcal W}_2:M_2(u,v)\leq 2C_0\varepsilon_0\}$. 
This is a closed subset of ${\mathcal W}_1$. 

As in the proof of Theorem \ref{GETheorem}, 
we employ the standard iteration method; 
namely, we set 
$(u_m,v_m)$ $(m=0,1,\dots)$ inductively 
by solving 
$\Box u=|v_{m-1}|^p$, 
$\Box v
=
|\partial_t u_{m-1}|^{(n+1)/(n-1)}$, 
with 
$(u(0),\partial_t u(0))=\varepsilon(f_1,g_1)$, 
$(v(0),\partial_t v(0))=\varepsilon(f_2,g_2)$. 
(Here, we have set $u_{-1}=v_{-1}\equiv 0$.) 
In the same way as in (\ref{202012161715}), 
it is possible to get, 
with a positive constant $C$ independent of $\varepsilon$, 
\begin{align}
\|&u_m(0)\|_2
+
\sum_{|\alpha|\leq 2}
\sum_{j=0}^n
\|
(\partial_j\Gamma^\alpha u_m)(0)
\|_2\\
&
+
\sum_{|\alpha|\leq 2}
\biggl(
\|
|D|^{\sigma(2p)}
\bigl(
(\Gamma^\alpha v_m)(0)
\bigr)
\|_2
+
\|
|D|^{\sigma(2p)-1}
\bigl(
(\partial_t\Gamma^\alpha v_m)(0)
\bigr)
\|_2
\biggr)
\leq
C\varepsilon
\nonumber
\end{align}
when $\varepsilon$ is sufficiently small. 
Moreover, 
owing to the assumption $p>1+3/(n-1)$, 
we obtain:
\begin{proposition}\label{202101221707}
For $m\geq 0$, the estimates
\begin{equation}\label{202101211828}
X_2(u_m)
\leq
C\varepsilon
+
C{\hat B}_2(v_{m-1})^p,
\quad
{\hat B}_2(v_m)
\leq
C\varepsilon
+
CX_2(u_{m-1})^{(n+1)/(n-1)}
\end{equation}
hold.
\end{proposition}
The proof of the latter inequality in (\ref{202101211828}) 
uses the argument in (\ref{202101071736}), 
(\ref{202101211831}), (\ref{202102241613}), 
(\ref{202101211833}), and (\ref{202102241616}) 
with $s=2p$. 
Moreover, suitably modifying the argument in 
(\ref{202101211838})--(\ref{202012211645}), 
we get the former one in (\ref{202101211828}). 

We see by Proposition \ref{202101221707} that 
there exist constants 
$C_0$ and $\varepsilon_0$ such that 
if $\varepsilon\leq\varepsilon_0$, 
then 
$(u_m,v_m)\in{\mathcal W}_2(2C_0\varepsilon_0)$, 
$m=0,1,\dots$ 
Moreover, 
repeating essentially the same argument as above, 
we can prove:
\begin{proposition}
For $m\geq 2$, the estimate
\begin{align}
M_1&(u_{m+1}-u_m,v_{m+1}-v_m)\\
&
\leq
C(2C_0\varepsilon_0)^{p-1}
{\hat B}_1(v_m-v_{m-1})
+
C(2C_0\varepsilon_0)^{2/(n-1)}
X_1(u_m-u_{m-1})\nonumber
\end{align}
holds.
\end{proposition}
The rest of the proof of Theorem \ref{202101221717} 
is obvious, and it is therefore omitted. 
We have finished the proof of Theorem \ref{202101221717}.
\section{}
This appendix is devoted to the proof of (\ref{202012231839}). 
In fact, we have only to repeat essentially the same argument 
as in \cite[Chapter 9]{Al2010} and \cite{HY2017}. 
The proof starts with the energy-type estimate 
for the inhomogeneous wave equation. 
Set $T_j=\partial_j+(x_j/|x|)\partial_t$, 
$j=1,2,3$. 
\begin{lemma}\label{202012241712}
Let $T>0$. 
For any $\eta>0$ there exists a constant $C=C_\eta$ 
such that smooth solutions $u(t,x)$ to $\Box u=F$ 
satisfy 
\begin{align}\label{202012241015}
\sup_{0<t<T}
\|&\partial u(t,\cdot)\|_2
+
\sum_{j=1}^3
\|
(1+|\tau-|x||)^{-1/2-\eta}T_j u
\|_{L^2((0,T)\times{\mathbb R}^3)}\\
&
\leq
C
\biggl(
\|\partial u(0,\cdot)\|_2
+
\int_0^T
\|F(\tau,\cdot)\|_2d\tau
\biggr).\nonumber
\end{align}
\end{lemma}
The proof of \cite[Lemma (energy inequality)]{Al2010}, 
which is based on the ghost weight technique, 
is obviously valid for that of (\ref{202012241015}). 
See \cite{Al2010} on page 92. 
It should be mentioned that a closely related estimate 
was obtained by Lindblad and Rodnianski 
in a different way \cite[Corollary 8.2]{LR2005}. 
See also \cite[(1.2)]{LNS}. 

It suffices to bound 
$\|\partial\Gamma^\alpha u_i^{(m)}(t)\|_2$, 
$\|\langle\tau-|x|\rangle^{-1/2-\eta}
T_j\Gamma^\alpha u_i^{(m)}\|_{L^2((0,T)\times{\mathbb R}^3)}$ 
$(|\alpha|\leq 3)$ only for $i=2$; the way of 
dealing with these norms for $i=1$ is similar 
(and a little simpler, in fact). 
In view of Lemma \ref{202012241712} 
we must deal with 
\begin{equation}
\int_0^t
\|
A_k^{ij,\alpha\beta}
(\partial_\alpha u_i^{(m-1)})
(\partial_\beta u_j^{(m-1)})
\|_{\Gamma,3,2}d\tau,\quad
k=1,2.
\end{equation}
We may focus on $k=2$; 
the other case is a little simpler to handle. 
In what follows, we again drop the superscript $(m-1)$ for simplicity. 
\subsection{$k=2$, $(i,j)=(1,1)$.} It is easy to get 
\begin{align}
\|&
A_2^{11,\alpha\beta}
(\partial_\alpha u_1)
(\partial_\beta u_1)
\|_{\Gamma,3,2}\\
&
\leq
\|
A_2^{11,\alpha\beta}
(\partial_\alpha u_1)
(\partial_\beta u_1)
\|_{\Gamma,3,2,\chi_1}
+
\|
A_2^{11,\alpha\beta}
(\partial_\alpha u_1)
(\partial_\beta u_1)
\|_{\Gamma,3,2,\chi_2}\nonumber\\
&
\leq
C(1+\tau)^{-1}
\|
\partial u_1(\tau)
\|_{\Gamma,3,2}^2
\leq
C(1+\tau)^{-1}X(u_1)^2,\nonumber
\end{align}
which yields
\begin{equation}
\int_0^t
\|
A_2^{11,\alpha\beta}
(\partial_\alpha u_1)
(\partial_\beta u_1)
\|_{\Gamma,3,2}
d\tau
\leq
C\bigl(\log(1+t)\bigr)X(u_1)^2.
\end{equation}
Here, we have used the Klainerman-Sobolev inequalities 
(\ref{202012211620})--(\ref{202012211648}) 
and 
the trace-type inequality (\ref{20201221933}) 
to deal with $\|\cdots\|_{\Gamma,3,2,\chi_1}$ 
and $\|\cdots\|_{\Gamma,3,2,\chi_2}$, 
respectively.  
\subsection{$k=2$, $(i,j)=(1,2)$.} 
In the same way as above, 
we get
\begin{equation}
\int_0^t
\|
A_2^{12,\alpha\beta}
(\partial_\alpha u_1)
(\partial_\beta u_2)
\|_{\Gamma,3,2}
d\tau
\leq
C(1+t)^\delta
X(u_1)Y(u_2).
\end{equation}
(As for this estimate, the null condition (\ref{Aass2}) plays no role.)
\subsection{$k=2$, $(i,j)=(2,2)$.} 
We deal with 
$\|
A_2^{22,\alpha\beta}
(\partial_\alpha u_2)
(\partial_\beta u_2)
\|_{\Gamma,3,2,\chi_i}$ 
for $i=1,2$, separately. 
For $i=1$, 
it is possible to get by the Klainerman-Sobolev inequalities 
(\ref{202012211620})--(\ref{202012211648}) 
\begin{equation}
\int_0^t
\|
A_2^{22,\alpha\beta}
(\partial_\alpha u_2)
(\partial_\beta u_2)
\|_{\Gamma,3,2,\chi_1}
d\tau
\leq
CY(u_2)^2.
\end{equation}
For $i=2$, the null condition (\ref{Aass3}) comes into play. 
We need:
\begin{lemma}\label{202012281821}
Suppose that 
$\{A^{\alpha\beta}:\alpha,\beta=0,1,2,3\}$ satisfies the null condition. 
Then, we have 
for smooth functions $w_i(t,x)$ $(i=1,2)$
\begin{equation}\label{202012241805}
|
A^{\alpha\beta}
(\partial_\alpha w_1)
(\partial_\beta w_2)
|
\leq
C
\biggl(
\sum_{j=1}^3
|T_j w_1|
\biggr)
|\partial w_2|
+
|\partial w_1|
\biggl(
\sum_{j=1}^3
|T_j w_2|
\biggr).
\end{equation}
\end{lemma}
For the proof, see \cite[pp.\,90--91]{Al2010}. 
It should be mentioned that 
a closely related inequality was obtained by 
Lindblad and Rodnianski \cite[(5.10)]{LR2005}. 
See also \cite[Lemma 2.3]{LNS}. 

We use (\ref{202012241805}) together with 
the useful idea of dyadic decomposition 
of the time interval $(0,\infty)$, 
and we also use (\ref{202012211743}) 
together with 
the trace-type inequality (\ref{20201221932}). 
The idea of how to use the decomposition of the time interval 
and the trace inequality in this way 
comes from \cite[p.\,363]{Sogge2003} and \cite[(3.38d)]{Sideris2000}, 
respectively. 
We then obtain 
\begin{align}
\int_0^\infty&
\|
A_2^{22,\alpha\beta}
(\partial_\alpha u_2)
(\partial_\beta u_2)
\|_{\Gamma,3,2,\chi_2}
d\tau\\
&
\leq
C
\bigl(
\sup_{t>0}
\langle t\rangle^{-\delta}
\|\partial u_2(t)\|_{\Gamma,3,2}
\bigr)\nonumber\\
&
\hspace{0.5cm}
\times
\biggl(
\sum_{\substack{|\alpha|\leq 3\\j=1,2,3}}
\sup_{t>0}\langle t\rangle^{-\delta}
\|\langle\tau-|x|\rangle^{-1/2-\eta}
T_j\Gamma^\alpha u_2\|_{L^2((0,t)\times{\mathbb R}^3)}
\biggr)\nonumber\\
&
\hspace{0.3cm}
+
\bigl(
\sup_{t>0}(1+t)^{-\delta}\|\partial u_2(t)\|_{\Gamma,3,2}
\bigr)^2\nonumber\\
&
\leq
CY(u_2)^2.\nonumber
\end{align}
See \cite[(3.61)--(3.66)]{{HY2017}} for details. 
The proof of (\ref{202012231839}) has been finished.
\section*{Acknowledgments} 
Special thanks go to the referee for a lot of 
valuable suggestions which 
have improved the presentation of this paper. 
It is also a pleasure to thank Professors 
Soichiro Katayama and Hideaki Sunagawa 
for their helpful comments. 
The first author was supported in part by 
the Grant-in-Aid for Scientific Research (C) (No.\,18K03365), 
Japan Society for the Promotion of Science (JSPS). 
\bibliographystyle{amsplain}

%
\end{document}